\documentclass{amsart}
\usepackage[utf8]{inputenc}  

\usepackage{amsfonts, amssymb}

\usepackage{amssymb,amsmath,amscd,euscript,verbatim,array}
\usepackage{geometry}
\usepackage{anysize}
\usepackage{fancyhdr}
\usepackage{indentfirst}
\usepackage{graphicx}
\usepackage{color}
\usepackage{ifpdf}
\usepackage{marginnote}

\newcounter{item}[section]

\newtheorem{Theorem}{Theorem}
\newtheorem{Definition}{Definition}[section]
\newtheorem{Proposition}[Definition]{Proposition}
\newtheorem{Corollary}[Definition]{Corollary}
\newtheorem{Lemma}[Definition]{Lemma}

\theoremstyle{definition}
\newtheorem{Remark}[Definition]{Remark}

\newcommand{\eps}{\varepsilon}

\newcommand{\bthm}{\begin{Theorem}}
	\newcommand{\ethm}{\end{Theorem}}
\newcommand{\bpr}{\begin{Proposition}}
	\newcommand{\epr}{\end{Proposition}}
\newcommand{\blm}{\begin{Lemma}}
	\newcommand{\elm}{\end{Lemma}}
\newcommand{\bex}{\begin{Exercise}}
	\newcommand{\eex}{\end{Exercise}}
\newcommand{\be}{\begin{equation}}
	\newcommand{\ee}{\end{equation}}
\newcommand{\beal}{\begin{aligned}}
	\newcommand{\enal}{\end{aligned}}
\newcommand{\brm}{\begin{Remark}}
	\newcommand{\erm}{\end{Remark}}

\begin{document}

\title[Moser's twist theorem revisited ]{Moser's twist theorem revisited }

\author{Yi Liu}
\address{School of Mathematics and Statistics, Beijing Institute of Technology, Beijing 100081, China}
\email{yiliu111@foxmail.com}

\author{Lin Wang}
\address{School of Mathematics and Statistics, Beijing Institute of Technology, Beijing 100081, China}
\email{lwang@bit.edu.cn}



\begin{abstract} {
Inspired by the work of Katznelson and Ornstein \cite{KO1}, we present a short way to achieve the almost optimal regularity in Moser's twist theorem. Specifically, for an integrable area-preserving twist map, the invariant circle with a given constant type frequency \( \alpha \) persists under a small perturbation (dependent on \( \alpha \)) of class \( C^{3+\eps} \). This result was initially established independently by Herman \cite{H1} and R\"{u}ssmann \cite{R2} in 1983. Our method differs essentially from their approaches.
}
\end{abstract}

\maketitle

\tableofcontents
\section{Introduction}
The existence of invariant circles (i.e., homotopically non-trivial invariant curves) for area-preserving twist maps was first established by Moser in his seminal work \cite{Mo1}, commonly referred to as \emph{Moser's twist theorem}. This groundbreaking result is also recognized as one of the foundational contributions to the KAM (Kolmogorov-Arnold-Moser) theory, where the map was initially required to be of class \( C^{333} \). The regularity condition was significantly relaxed to \( C^{4+\varepsilon} \) by R\"{u}ssmann \cite{R1}, where $\varepsilon$ denotes, once and for all, an arbitrarily small positive constant. Subsequently, Moser \cite[Page 53]{Mo2} asserted that \( C^{3+\varepsilon} \) regularity suffices for the existence of invariant circles by employing more refined estimates within the KAM framework. The full proof of this assertion was independently provided by R\"{u}ssmann \cite{R2} and Herman \cite{H1}, each utilizing distinct methodologies. In \cite{R2}, R\"{u}ssmann introduced a novel iteration process that yields invariant circles which are merely continuous. In contrast, Herman's approach in \cite{H1} leverages the Schauder fixed point theorem, building upon his celebrated work on circle diffeomorphisms (\cite{H11}). Although Herman's method is restricted to invariant circles with constant type frequency, it achieves superior smoothness for these circles, specifically \( C^{2+\varepsilon} \).

Herman \cite{H1} also constructed an example demonstrating that each invariant circle can be destroyed by a small $C^{3-\varepsilon}$ perturbation. This implies that $C^{3+\varepsilon}$ regularity is almost optimal. Furthermore, Herman \cite{H33}, using a result by Meyer \cite{Me}, achieved the optimal regularity result. Specifically, he showed that an invariant circle with a constant type frequency persists under small perturbations in the $C^{3}$ topology.

Inspired by the work of Katznelson and Ornstein \cite{KO1}, we propose a novel approach to determine the almost optimal regularity in Moser's twist theorem. This method diverges from both R\"{u}ssmann's and Herman's techniques. Similar to Herman \cite{H1}, our goal is to achieve the smoothness ($C^2$ plus the H\"{o}lder condition) for invariant circles with constant type frequency. However, this approach alone appears insufficient to reduce the regularity requirement from \( C^{3+\varepsilon} \) to the optimal \( C^{3} \). We plan to address this limitation in future research.

Let \( \mathbb{T} := \mathbb{R}/\mathbb{Z} \) denote the flat circle. Let \( |\cdot| \) represent the Euclidean metric in \( \mathbb{R}^d \). The flat metric on \( \mathbb{T} \) is defined as
\begin{equation}\label{flamet}
\|\cdot\| := \inf_{p \in \mathbb{Z}} |\cdot + p|.
\end{equation}
Let \( \alpha \) be an irrational number, expressed in its continued fraction form as \( \alpha = [a_0, a_1, \ldots] \). The convergents \( \left\{ \frac{p_i}{q_i} \right\}_{i \in \mathbb{N}} \) of \( \alpha \) are given by \( \frac{p_i}{q_i} = [a_0, \ldots, a_i] \). Without loss of generality, we assume \( \alpha \in (0, 1) \) and \( q_n > 0 \) for all \( n \in \mathbb{N} \). An irrational number \( \alpha \) is said to be of constant type if all coefficients \( a_i \) in its continued fraction are bounded. We denote this bound by \( A_\alpha \). Let \( \mathcal{C} \) represent the set of constant type irrational numbers in \( (0, 1) \). It is well-known that \( \mathcal{C} \) has zero Lebesgue measure but full Hausdorff dimension.

Given $k \in [0, +\infty) \cup \{+\infty\}$, let $C^k(\mathbb{T})$ denote the space of $C^k$ functions on the circle $\mathbb{T}$. Given any $u\in C^k(\mathbb{T})$. Let us recall the $C^k$-norm:
\[\|u\|_{C^k}:=\max_{|r|\leq k}\max_{x\in \mathbb{T}}|D^{r}u(x)|,\]
where $r:=(r_1,\ldots,r_n)$, $|r|:=r_1+\cdots+r_n$.

We consider a special class of generating functions \( \{G_n\}_{n \in \mathbb{N}} \) defined by
\begin{equation}\label{epss}
G_n(x, x') := \frac{1}{2}(x - x')^2 + \frac{1}{q_n^{4+\varepsilon}} V(q_n x'),
\end{equation}
where  \( V \in C^{4+\varepsilon}(\mathbb{T}) \) and \( q_n \) denotes the denominator of the \(n\)-th convergent of \( \alpha \). This means that the construction of \( G_n \) depends on \( \alpha \); for notational simplicity, we suppress the subscript \( \alpha \) in our notation. Given $n\in \mathbb{N}$, let \( F_n: \mathbb{T} \times \mathbb{R} \to \mathbb{T} \times \mathbb{R} \) be the exact area-preserving twist map generated by \( G_n \).

\begin{Theorem}\label{mathe}
For each \( \alpha \in \mathcal{C} \), there exists a sufficiently large \( n \in \mathbb{N} \) such that \( F_n \) admits an invariant circle with frequency \( \alpha \) and the preserved invariant circle can be represented as a graph of the function of class \( C^{2+\varepsilon'} \) for any \( \varepsilon' < \varepsilon \).
\end{Theorem}

It is straightforward to verify that Theorem \ref{mathe} remains valid if the perturbation term in \( G_n \) is replaced by \( U_n \in C^{4+\varepsilon}(q_n \mathbb{R}/\mathbb{Z}) \), provided \( \|U_n\|_{C^{4+\varepsilon}} \) is bounded independently of \( n \). It is also noticed that  $\varepsilon'=\varepsilon$ is still true (\cite{H1}).

Let \( T(x, y) = (x + y, y) \). By definition, for each \( \varepsilon' < \varepsilon \), we have
\[
\|F_n - T\|_{C^{3+\varepsilon'}} \to 0 \quad \text{as} \quad n \to \infty.
\]
According to \cite{H1}, any invariant circle can be destroyed by a \( C^{3-\varepsilon} \) perturbation. Thus, Theorem \ref{mathe} is nearly optimal (up to \( \varepsilon \)).
The proof of Theorem \ref{mathe} combines the circle diffeomorphism technique developed by Katznelson and Ornstein \cite{KO1} with the action-minimizing method introduced by Mather \cite{FoM}.

Moser \cite{Mo4} observed that every $C^\infty$ twist diffeomorphism corresponds to the time-one map of a suitable $C^\infty$ periodic Hamiltonian system, although the converse does not hold. To ensure the validity of the converse, the condition that there are no conjugate points must also be imposed. Consequently, a result analogous to Theorem \ref{mathe} can be directly derived for Hamiltonian systems with two degrees of freedom. More precisely, let $\mathbb{T}^2$ be the 2-dimensional flat torus. We denote by $\mathcal{D} \subset \mathbb{R}^2$ the set of frequency vectors $\omega \in \mathbb{R}^2$ satisfying the constant-type Diophantine condition: there exists a constant $D > 0$ such that for all non-zero integer vectors $k \in \mathbb{Z}^2 \setminus \{0\}$,
\begin{equation}\label{eq:diophantine}
	|\langle \omega, k \rangle| \geq \frac{D}{|k|}.
\end{equation}

Consider the integrable Hamiltonian $H_0: \mathbb{T}^2 \times \mathbb{R}^2 \to \mathbb{R}$ given by
\[
H_0(q, p) = \frac{1}{2} |p|^2, \quad p \in \mathbb{R}^2.
\]
\begin{Corollary}\label{macor1}
	For each $\omega \in \mathcal{D}$, there exists a sequence of perturbed Hamiltonians $\{H_n\}_{n\in\mathbb{N}}$ of class $C^{4+\varepsilon}$  such that for any $\varepsilon' < \varepsilon$:
	\begin{enumerate}
		\item $\|H_n - H_0\|_{C^{4+\varepsilon'}} \to 0$ as $n \to +\infty$;
		\item each $H_n$ admits an invariant torus with frequency $\omega$;
		\item the preserved invariant torus is a graph of class $C^{2+\varepsilon'}$.
	\end{enumerate}
\end{Corollary}

While Corollary \ref{macor1} is not a new result,  it has appeared in many literatures, either explicitly or implicitly (see \cite{A,Sa,TL} for instance). We include it here as a corollary of Theorem \ref{mathe} merely to illustrate that this result can be obtained through a novel approach distinct from existing literatures. Under the consideration that an individual KAM torus with constant frequency persists, the current findings on the minimal smoothness requirement ($C^{2d}$ plus the Dini condition) for the Hamiltonian function are presented in \cite{A,TL}.  By \cite{Bou, Kou}, it is also possible, through a modified KAM approach, to identify nearly integrable Hamiltonian systems that admit \( d \)-dimensional invariant tori with positive Lebesgue measure under \( C^{2d+\varepsilon} \) perturbations. Conversely, as shown in \cite{CW,W2}, any \( d \)-dimensional Lagrangian invariant torus can be destroyed by a small \( C^{2d-\varepsilon} \) perturbation.

\vspace{1em}

\noindent\textbf{Strategy of the proof of Theorem \ref{mathe}.} Katznelson and Ornstein \cite{KO1} introduced a criterion for the existence of invariant circles (referred to as Criterion 1 in this paper), expressed in terms of certain vectors with endpoints on orbit segments (called Type-I chords here; see Definition~\ref{kvec1}). However, working directly with such chords makes it difficult to achieve almost optimal result on the regularity of the perturbation for twist maps. We therefore introduce a more refined class of chords (called Type-II chords; see Definition~\ref{kvec3}), whose definition imposes an additional constraint on the ``difference of positions'' of their endpoints on the minimal configuration.

Based on these Type-II chords, we formulate a new criterion for invariant circles (Criterion 2) and prove its equivalence to Criterion 1 (Lemma~\ref{3lam}). To connect the construction with circle diffeomorphism theory, we further reduce Criterion 2 to Criterion 3, which is a boundedness condition on \(\widetilde{\mathbf{K}}^0\) (defined  in (\ref{bigK})) related to a certain multiplicative cocycle. This reduction is carried out in Section 3.

In light of the idea to prove the Denjoy inequality for circle diffeomorphisms, Katznelson and Ornstein \cite{KO1} established a corresponding result for twist maps  (Proposition~\ref{ko11}), which implies the quantity \(\widetilde{\mathbf{K}}^0\) can be controlled by a certain distortion \(\widetilde{\mathbf{K}}^1\) (defined in (\ref{distK})). Introducing the second and third-order difference quotients \(\nabla^1\) and \(\nabla^2\) of the minimal configuration (see definitions (\ref{na11}) and (\ref{na22})), and under suitable hypotheses on \(\nabla^1,\nabla^2\) (conditions (\ref{f1x}) and (\ref{f2x})), one obtains exponential decay of \(\widetilde{\mathbf{K}}^1\) (Lemma~\ref{last}).  Boundedness of \(\widetilde{\mathbf{K}}^0\) then follows from Proposition~\ref{ko11}.

The proof of Lemma~\ref{last} is most delicate part in this paper. Establishing the exponential decay of \(\widetilde{\mathbf{K}}^1\) requires several steps. We first prove a monotonicity property for \(\nabla^1\) (Lemma~\ref{A1}). When the orbit segment is short relative to the scale of the chord-pair involved in \(\nabla^1\), we show exponential decay of \(\widetilde{\mathbf{K}}^1\) (Lemma~\ref{keyy1}). Using the assumption on \(\nabla^2\), we then extend this decay to the case where the length of the orbit segment matches the scale involved in  \(\nabla^1\) (Lemma~\ref{karrr1}). Finally we obtain, a posteriori, a sharpened exponential decay estimate for \(\widetilde{\mathbf{K}}^1\) (Proposition~\ref{L1}). These arguments are contained in Section 4.

The last step consists of an induction argument showing that there exists a scale \(\kappa_0\) such that, on all scales finer than \(\kappa_0\), the conditions on \(\nabla^1\) and \(\nabla^2\) ((\ref{f1x}) and (\ref{f2x})) are automatically satisfied (Lemma~\ref{las1}). The initial scale \(\kappa_0\) depends in a precise way on the size of the perturbation of the system (Lemma~\ref{initi}). This is presented in Section~5, together with a basic estimate for the vertical-coordinate difference between consecutive points of the orbit (Lemma~\ref{disne}) given in Section~2. To deduce the final result from the induction hypothesis, we need to control the growth of \(\nabla^1\) in relation to \(\nabla^2\) (Lemma~\ref{SSH}) as the scale of the chord increases, and then combine this with the refined estimate for \(\nabla^2\) (see (\ref{keyfor})) obtained by Katznelson and Ornstein (\cite[Corollary 3.1]{KO1}) to complete the proof of Theorem~\ref{mathe}.

Roughly speaking, the main idea to prove Theorem~\ref{mathe}  is to prove that an Aubry-Mather set which ``looks like'' an invariant circle must indeed be an invariant circle. The precise meaning of ``looks like'' is encoded by the second and third-order difference quotients \(\nabla^1\) and \(\nabla^2\) of the minimal configuration. Through a bootstrap-type procedure, such an Aubry-Mather set cannot contain any ``gaps'' and therefore necessarily satisfies Criterion 1.

\vspace{2em}

\noindent\textbf{Acknowledgement.}   The authors would like to thank Zhicheng Tong for his helpful comments.  This work is partially under the support of National Natural Science Foundation of China (Grant No. 12122109) and the Beijing Natural Science Foundation (Grant No. QY25241).

\vspace{2em}

\section{The minimal configuration with frequency $\alpha$}\label{S2}
Let us recall some basic notions as preparations. Let $\tilde{F}$ be a diffeomorphism of $\mathbb{R}^2$ denoted by $\tilde{F}(x,y)=(X(x,y),Y(x,y))$. Let $\tilde{F}$ satisfy:
\begin{itemize}
	\item {\it Lift condition:} $\tilde{F}$ is isotopic to the identity;
	\item {\it Twist condition:} the map $\psi:(x,y)\mapsto(x,X(x,y))$ is a diffeomorphism of $\mathbb{R}^2$;
	\item {\it Exact symplectic:} there exists a real valued function $H$ on $\mathbb{R}^2$ with $H(x+1,X+1)=H(x,X)$ such that
	\[YdX-ydx=dH.\]
\end{itemize}
Then $\tilde{F}$ induces a map on the cylinder denoted by $f$: $\mathbb{T}\times\mathbb{R}\mapsto \mathbb{T}\times\mathbb{R}$. $f$ is called an exact
area-preserving  twist map. The function  $H$: $\mathbb{R}^2\rightarrow\mathbb{R}$ is called
a generating function of $\tilde{F}$, namely $\tilde{F}$
is generated by the following equations
\begin{equation*}
	\begin{cases}
		y=-\partial_1 H(x,X),\\
		Y=\partial_2 H(x,X),
	\end{cases}
\end{equation*}
where $\tilde{F}(x,y)=(X,Y)$.

The function $\tilde{F}$ gives rise to a dynamical
system whose orbits are given by the images of points of $\mathbb{R}^2$
under the successive iterates of $\tilde{F}$. The orbit of the point
$(x_0,y_0)$ is the bi-infinite sequence
\[\{...,(x_{-k},y_{-k}),...,(x_{-1},y_{-1}),(x_0,y_0),(x_1,y_1),...,(x_k,y_k),...\},\]
where $(x_k,y_k)=\tilde{F}(x_{k-1},y_{k-1})$. The sequence
\[(...,x_{-k},...,x_{-1},x_0,x_1,...,x_k,...)\] denoted by ${\bf{x}}:=(x_i)_{i\in\mathbb{Z}}$ is called a
stationary configuration associated to $\tilde{F}$ if it satisfies the discrete Euler-Lagrange equation
\[\partial_1 H(x_i,x_{i+1})+\partial_2 H(x_{i-1},x_i)=0,\ \text{for\ every\ }i\in\mathbb{Z}.\]
Given a sequence of points $(z_i,...,z_j)$, we can associate its
action
\[W(z_i,...,z_j)=\sum_{i\leq s<j}H(z_s,z_{s+1}).\] A configuration $(x_i)_{i\in\mathbb{Z}}$
is called minimal if for any $i<j\in \mathbb{Z}$, the segment
$(x_i,...,x_j)$ minimizes $W(z_i,...,z_j)$ among all segments
$(z_i,...,z_j)$ of the configuration  satisfying $z_i=x_i$ and
$z_j=x_j$. Then every minimal configuration is a
stationary configuration. For every minimal configuration $\bold{x}=(x_i)_{i\in\mathbb{Z}}$, the limit
\[\rho(\bold{x})=\lim_{n\rightarrow\infty}\frac{x_{i+n}-x_i}{n}\]
exists and it is independent of $i\in\mathbb{Z}$. $\rho(\bold{x})$ is called
the frequency of $\bold{x}$. We recall the orientation-preserving property as follows (see \cite[Proposition 11.2.4]{KH}).
\begin{Proposition}\label{OP}
If $\rho(\bold{x})$ is irrational, then for $n_1,n_2,m_1,m_2\in \mathbb{Z}$,
\[n_1\rho(\bold{x})+m_1<n_2\rho(\bold{x})+m_2\]
if and only if
\[x_{n_1}+m_1<x_{n_2}+m_2.\]
\end{Proposition}

For simplicity, we  use, once and for all,  $u\lesssim v$ (resp.
$u\gtrsim v$) to denote $u\leq Cv$ (resp. $u\geq Cv$), and $u\sim v$ to denote $\frac{1}{C}v\leq u\leq Cv$ for some
positive constant $C$.

Based on the construction of $G_n$, for each  $i\in \mathbb{Z}$, we have
\[x_{i+1}-x_{i}=y_i,\quad y_{i+1}-y_i=\frac{1}{q_n^{3+\varepsilon}}V'(q_n(x_i+y_i)).\]It is clear to see
\begin{Lemma}\label{disne}
Let $(x_i)_{i\in\mathbb{Z}}$ be a minimal configuration with rotation number $\alpha$.
Then \[|y_{i+1}-y_i|\lesssim \frac{1}{q_n^{3+\varepsilon}},\quad \forall \ i\in\mathbb{Z}.\]
\end{Lemma}

\vspace{2em}

\section{Criterions on the existence of the invariant circle }\label{S3}

Let \( F:\mathbb{T}\times\mathbb{R}\to\mathbb{T}\times\mathbb{R} \) be  an exact area-preserving twist map defined on \( \mathbb{T}\times\mathbb{R} \).  Take an irrational number \( \alpha \). It is known (see \cite{B,FoM}) that there exists an Aubry--Mather set \( \widetilde{\mathcal{M}}_{\alpha}\subset \mathbb{T}\times\mathbb{R} \).

Let $\pi:\mathbb{T}\times\mathbb{R}\to\mathbb{T}$ be the standard projection. Briefly, the projected set $\pi \widetilde{\mathcal{M}}$ onto \( \mathbb{T} \)  consists of all minimal configurations with rotation number \( \alpha \). For simplicity, we do not distinguish here between a minimal configuration on \(\mathbb{T}\) and its lift to \(\mathbb{R}\). Denote by \( H:\mathbb{R}^2\to\mathbb{R} \) the generating function of the lift \( \tilde{F} \).
More precisely, let \( \mathbf{x}=(x_i)_{i\in\mathbb{Z}} \) be a minimal configuration generated by \( H \). Define
\[
y_0:=-\partial_1 H(x_0,x_1).
\]
Let
\[
\mathcal{O}:= \{ F^{i}(x_0,y_0)\}_{i\in\mathbb{Z}},
\]
where \( \mathcal{O} \) is a complete orbit passing through \( (x_0,y_0)\).
If the rotation number \( \rho(\mathbf{x}) \) is irrational, and we set \( \alpha=\rho(\mathbf{x}) \), then \( \widetilde{\mathcal{M}}_{\alpha} \) can be expressed as
\[
\widetilde{\mathcal{M}}_{\alpha}= \overline{\mathcal{O}},
\]
where \( \overline{\mathcal{O}} \) denotes the closure of \( \mathcal{O} \) in the topology induced by the Riemannian metric on \( \mathbb{T}\times\mathbb{R} \).

The mapping \( x_i \mapsto i\alpha \) is orientation-preserving from \( \mathcal{O} \) to the circle \( \mathbb{T} \). A vector \( v \) is called a chord if its endpoints lie on some orbit segment of \( \mathcal{O} \). Denote by \( \lambda_1(v) \) the length of the arc between the corresponding multiples of \( \alpha \) on \( \mathbb{T} \). Concretely, if the chord \( v \) connects \( z_i=(x_i,y_i) \) and \( z_j=(x_j,y_j) \), then
\[
\lambda_1(v)=\|(j-i)\alpha\|,
\]
where \( \|\cdot\| \) is defined by (\ref{flamet}). It is clear that \( \lambda \) is invariant under \( F \). For clarity, \( \lambda_1(v) \) is called the \(\lambda\)-length of $v$ on \( \mathcal{O} \).
Under the same assumption that the chord \( v \) connects \( z_i=(x_i,y_i) \) and \( z_j=(x_j,y_j) \), we also write
\[
\theta(v):=x_j-x_i, \qquad \Theta_1(v):=\|\theta(v)\|.
\]
The quantity \( \Theta_1(v) \) is called the \(\vartheta\)-length of $v$ on \( \mathcal{O} \).

\begin{Remark}
	By the definitions of the $\lambda$-length and the $\vartheta$-length of $v$ on $\mathcal{O}$, for two chords $v_1$ and $v_2$ on $\mathcal{O}$, even if their $\lambda$-distances coincide, their $\vartheta$-lengths may still differ.
	In $\widetilde{\mathcal{M}}_\alpha$ we fix, once and for all, an arbitrary complete orbit $\mathcal{O}$.
	All subsequent discussions about a chord $v$ refer to the vectors that arise from some orbit segment whose endpoints belong to this chosen complete orbit $\mathcal{O}$.
\end{Remark}

Based on the orientation-preserving property (see Proposition \ref{OP}),  $\lambda_1(v)$ and $\Theta_1(v)$ can be related with each other in the sense of average. Recall \( q_n \) denotes the denominator of the \(n\)-th convergent of \( \alpha \).

\begin{Proposition}\label{eveco}Let $v$ be a chord. If $\lambda_1(v)\geq\frac{2}{q_N}$ for $N$ large enough, then
\[\left|\frac{1}{q_N}\sum_{j=0}^{q_N-1}\Theta_1(F^j(v))-\lambda_1(v)\right|\leq \frac{1}{q_N}.\]
\end{Proposition}

\subsection{$\kappa$-chords}
Since $\alpha$ is irrational, then $\{\|m\alpha\|\}_{m\in \mathbb{Z}}$ is dense in $\mathbb{T}$. It means that for each $\kappa\in \mathbb{N}$, there exists $m:=m(\kappa)$ such that $\|m\alpha\|\in  [\frac{1}{2^{\kappa+4}},\frac{1}{2^\kappa}]$.
Given $n_0\in \mathbb{N}$, one can find $\kappa_0:=\kappa_0(n_0)$ such that
\begin{equation}\label{qn}
\|q_{n_0}\alpha\|\in \left[\frac{1}{2^{\kappa_0+4}},\frac{1}{2^{\kappa_0}}\right].
\end{equation}
In fact, by Dirichlet's approximation,
\begin{equation}\label{DA}
\frac{1}{(a_{n+1}+2)q_n}<\|q_n\alpha\|<\frac{1}{a_{n+1}q_n}.
\end{equation}
We choose
\begin{equation}\label{k0}
\kappa_0=\mathrm{log}_2a_{n_0+1}+\mathrm{log}_2q_{n_0}.
\end{equation}
Then we get
\[\frac{1}{2^{\kappa_0+4}}\leq \frac{1}{(a_{n_0+1}+2)q_{n_0}}<\frac{1}{a_{n_0+1}q_{n_0}}\leq \frac{1}{2^{\kappa_0}}.\]

Define $\varphi(m):=\mathrm{log}_2a_{m+1}+\mathrm{log}_2q_{m}$. Note that $q_{m+1}=a_{m+1}q_m+q_{m-1}$. A direct calculation shows
\begin{equation}
0<\varphi(m+1)-\varphi(m)\leq  \mathrm{log}_2A_\alpha+1.
\end{equation}
We also define
\begin{equation}\label{phiex}
\phi(m):=[\varphi(m)].
\end{equation}
We use, once and for all,  $[c]$ to denote the largest integer not greater than $c$ for $c\in \mathbb{R}$. It is clear to see that $\phi:\mathbb{N}\to\mathbb{N}$ is nondecreasing. In general, $\phi$ may not be surjective. Meanwhile, $\phi:\mathbb{N}\to\mathbb{N}$ may not be strict increasing. Consequently, $\phi^{-1}:\mathbb{N}\to \mathbb{N}$ may not be well-defined. Denote $S:=\phi(\mathbb{N})$. We introduce $\psi$ as a substitution of $\phi^{-1}:\mathbb{N}\to \mathbb{N}$. For each $\kappa\notin S$, there exists a maximal $\hat{\kappa}\in S$ satisfying $\hat{\kappa}<\kappa$.  Note that
\begin{equation}\label{disp3}
\phi(m+3)>\phi(m).
\end{equation}
We have
\[\kappa-2\leq \hat{\kappa}<\kappa.\]Denote
\[\check{\kappa}:=\left\{\begin{array}{ll}
	\hspace{-0.4em} \kappa,\quad \kappa\in S,\\
	\hspace{-0.4em} \hat{\kappa},\quad \kappa\notin S.
\end{array}\right.\]
The function $\psi:\mathbb{N}\to \mathbb{N}$ is defined as
\begin{equation}\label{phius}
\psi(\kappa):=\min\{m\in \mathbb{N}\ |\ \phi(m)=\check{\kappa}\}.
\end{equation}
It follows that $\kappa\to+\infty$ if and only if $\psi(\kappa)\to+\infty$. In the following, we define two types of $\kappa$-chords. The first one was introduced by \cite{KO1}.

\begin{Definition}[Type-I]\label{kvec1}
	We call $v$ a Type-I $\kappa$-chord, denoted by $v\in (\kappa)_{\text{I}}$, if  $\lambda_1(v)\in [\frac{1}{2^{\kappa+4}},\frac{1}{2^\kappa}]$.
\end{Definition}
To achieve the almost optimal regularity in Theorem \ref{mathe}, we need another kind of $\kappa$-chord.
\begin{Definition}[Type-II]\label{kvec3}
	Assume $v$ connects $(x_i,y_i)$ and  $(x_j,y_j)$. We call $v$ a Type-II $\kappa$-chord, denoted by $v\in (\kappa)_{II}$, if $v\in (\kappa)_{I}$ and
	\[{\frac{1}{(A_{\alpha}+1)^6}q_{n_\kappa}\leq|i-j|\leq (A_{\alpha}+1)^6q_{n_\kappa},\quad \lambda_1(v)\in [\|q_{n_\kappa}\alpha\|,2^4\|q_{n_\kappa}\alpha\|]}\]  where $n_{\kappa}:=\psi(\kappa)$.
\end{Definition}
For the  convergents $\{\frac{p_n}{q_n}\}_{n\in \mathbb{N}}$ of $\alpha$, there holds for $n\in\mathbb{N}$,
\begin{equation}\label{D2}
q_{n+2}> 2q_n,\quad (A_\alpha+1)^n>q_n> (\sqrt{2})^n.
\end{equation}
Given $\kappa\in \mathbb{N}$ large enough,
we denote, once and for all,
\begin{equation}\label{Z11}
\gamma_0:=\left[\mathrm{log}_22(A_\alpha+2)\right]+4,\quad \bar{N}(\kappa):=n_\kappa+2\gamma_0,
\end{equation}
\begin{equation}\label{Z12}
\tilde{N}(\kappa):=n_\kappa-2\gamma_0,
\end{equation}
where $n_\kappa=\psi(\kappa)$.
By (\ref{DA}), (\ref{D2}), (\ref{Z11}) and (\ref{Z12}), a direct calculation shows that if $v\in (\kappa)_{II}$, then
\begin{equation}\label{Z13}
\frac{2}{q_{\bar{N}(\kappa)}}<\lambda_1(v)<\frac{1}{4q_{\tilde{N}(\kappa)}}.
\end{equation}

\begin{Remark}\label{RGX}
	\begin{itemize}
		\
		\item It is easy to verify that the set $(\kappa)_{II}$ is nonempty. In fact, the chord $v$ joining $(x_i, y_i)$ and $(x_{i+q_{n_{\kappa}}}, y_{i+q_{n_{\kappa}}})$ belongs to $(\kappa)_{II}$.
		
		\item We illustrate the notations with two explicit constant-type irrational numbers.
		
		\begin{itemize}
			\item For $\alpha$ equal to the golden ratio \(\displaystyle \frac{\sqrt{5}-1}{2}=[0;1,1,\dots]\), we have \(A_\alpha=1\) and \(\gamma_{0}= \bigl[\log_2\bigl(2(A_\alpha+2)\bigr)\bigr]+4 = 6\). Because \(0<\varphi(m+1)-\varphi(m)\le 1\), the map \(\phi\) is surjective: \(\phi(\mathbb N)=\mathbb N\). Simple computations give \(\psi(1)=1,\;\psi(5)=9,\;\psi(10)=17,\;\psi(100)=146,\;\psi(500)=722\).
			
			\item For $\alpha=\sqrt{2}-1=[0;2,2,\dots]$, we have \(A_\alpha=2\) and \(\gamma_{0}= \bigl[\log_2\bigl(2(A_\alpha+2)\bigr)\bigr]+4 = 7\). Here \(\phi\) is not surjective; indeed \(\phi(\mathbb N)=\mathbb N\setminus \mathcal{T}\) with \(\mathcal{T}=\{4,9,14,18,23,28,\dots\}\), where the gap between two consecutive elements of \(\mathcal{T}\) is either \(4\) or \(5\). Calculations yield \(\psi(1)=1,\;\psi(4)=\psi(3)=4,\;\psi(28)=23,\;\psi(100)=80,\;\psi(500)=395\).
		\end{itemize}
	\end{itemize}
\end{Remark}

\subsection{Criteria}

We consider the following quantities
\[\Lambda_\iota(\kappa):=\max\left\{\sup_{v\in (\kappa)_\iota}\frac{\lambda_1(v)}{\Theta_1(v)},\sup_{v\in (\kappa)_\iota}\frac{\Theta_1(v)}{\lambda_1(v)}\right\},\quad \iota\in \{I,II\}.\]
Roughly speaking, $\Lambda_I(\kappa)$ can be viewed as an analog of the derivative of the conjugate function in the theory of circle diffeomorphisms. Based on this analogy, we derive a criterion for the existence of an invariant circle with a prescribed frequency. As demonstrated by Katznelson and Ornstein \cite{KO1}, for a given $n \in \mathbb{N}$, the map $F_n$ generated by $G_n$ admits an invariant circle with rotation number $\alpha$ if and only if the following condition is satisfied.
\begin{itemize}
	\item {\bf Criterion 1:}    The sequence $\{\Lambda_I(\kappa)\}_{\kappa\in \mathbb{N}}$ is  bounded.
\end{itemize}
Obviously, we only need to consider the case with $\kappa$ large enough. For $\kappa$ large enough,  if $v=v'+v''$, then we have $\lambda_1(v)=\lambda_1(v')+\lambda_1(v'')$ and $\Theta_1(v)=\Theta_1(v')+\Theta_1(v'')$. It follows that $\frac{\lambda_1(v)}{\Theta_1(v)}$ (resp. $\frac{\Theta_1(v)}{\lambda_1(v)}$) is a convex combination of $\frac{\lambda_1(v')}{\Theta_1(v')}$ and $\frac{\lambda_1(v'')}{\Theta_1(v'')}$ (resp. $\frac{\Theta_1(v')}{\lambda_1(v')}$ and $\frac{\Theta_1(v'')}{\lambda_1(v'')}$). Hence, $\Lambda_I(\kappa)$ is nondecreasing (up to a subsequence). It suffices to verify the boundedness of $\{\Lambda_{I}(\kappa)\}_{\kappa\in S}$. Moreover, one can show that $\{\Lambda_I(\kappa)\}_{\kappa\in S}$ is even controlled by $\{\Lambda_{II}(\kappa)\}_{\kappa\in \mathbb{N}}$, for which one needs to show the following
\begin{Lemma}\label{3lam}
For each \( v \in (\kappa)_{I} \) with \( \kappa \in S \), there exist \( v' \in (\kappa')_{II} \) and \( v'' \in (\kappa'')_{II} \) such that
\[
\frac{1}{2^4 (A_\alpha + 1)^8} \frac{\lambda_1(v'')}{\Theta_1(v'')}
< \frac{\lambda_1(v)}{\Theta_1(v)}
< (A_\alpha + 1)^8  \frac{\lambda_1(v')}{\Theta_1(v')},
\]
where \( A_\alpha \) denotes an upper bound for all the partial quotients \( a_i \) in the continued fraction expansion of \( \alpha \).
\end{Lemma}

\begin{proof}
Fix $n\in \mathbb{N}$. There holds (see for instance \cite[V. (7.6)]{H11})
\[2<\frac{\|q_n\alpha\|}{\|q_{n+2}\alpha\|}<(a_{n+2}+1)(a_{n+3}+1).\]
It follows that for $j\in \mathbb{N}$,
\[2^j\|q_{n+2j}\alpha\|<\|q_n\alpha\|<(A_\alpha+1)^{2j}\|q_{n+2j}\alpha\|.\]
In particular, we have
\[\|q_{n+8}\alpha\|<\frac{1}{2^4}\|q_n\alpha\|<\frac{1}{2^{\kappa+4}}<\frac{1}{2^\kappa}<2^4\|q_n\alpha\|<\|q_{n-8}\alpha\|<(A_\alpha+1)^{8}\|q_n\alpha\|.\]
Note that for  $\kappa=\phi(n)\in S$, $\|q_n\alpha\|\in [\frac{1}{2^{\kappa+4}},\frac{1}{2^\kappa}]$.
\begin{itemize}
	\item {\it Proof of $\frac{\lambda_1(v)}{\Theta_1(v)}< (A_\alpha + 1)^8 \frac{\lambda_1(v')}{\Theta_1(v')}$.}
\end{itemize}
Given $v\in (\kappa)_{I}$ with $\kappa=\phi(n)$, if $\lambda_1(v)\in [\|q_n\alpha\|, \frac{1}{2^\kappa}]$, then we only need to choose $v'\in (\kappa)_{II}$ satisfying $\lambda_1(v')=\|q_n\alpha\|$ such that
\[\frac{\lambda_1(v)}{\Theta_1(v)}< 2^4\frac{\|q_n\alpha\|}{\Theta_1(v')}= 2^4\frac{\lambda_1(v')}{\Theta_1(v')}.\]
In fact, by the orientation-preserving property (see Proposition \ref{OP}), $\Theta_1(v')<\Theta_1(v)$.

Otherwise, if  $\lambda_1(v)\in [\frac{1}{2^{\kappa+4}}, \|q_n\alpha\|)$, one choose  $v'\in (\kappa')_{II}$ with $\kappa'=\phi(n+8)$ satisfying $\lambda_1(v')=\|q_{n+8}\alpha\|$. Then we still have $\Theta_1(v')<\Theta_1(v)$ from the orientation-preserving property. Moreover,
\[\frac{\lambda_1(v)}{\Theta_1(v)}\leq \frac{\|q_n\alpha\|}{\Theta_1(v')}<(A_\alpha+1)^8\frac{\|q_{n+8}\alpha\|}{\Theta_1(v')}\leq(A_\alpha+1)^8\frac{\lambda_1(v')}{\Theta_1(v')}.\]
\begin{itemize}
	\item {\it Proof of $\frac{1}{2^4 (A_\alpha + 1)^8} \frac{\lambda_1(v'')}{\Theta_1(v'')}< \frac{\lambda_1(v)}{\Theta_1(v)}$.}
\end{itemize}		
Given $v\in (\kappa)_{I}$ with $\kappa=\phi(n)$,  we only need to choose $v''\in (\kappa'')_{II}$ with  $\kappa''=\phi(n-8)$ satisfying $\lambda_1(v'')=\|q_{n-8}\alpha\|$. Similarly, by the orientation-preserving property (see Proposition \ref{OP}), $\Theta_1(v'')>\Theta_1(v)$. It follows that
\[\frac{1}{2^4(A_\alpha+1)^{8}}\frac{\lambda_1(v'')}{\Theta_1(v'')}<\frac{1}{2^4}\frac{\|q_n\alpha\|}{\Theta_1(v'')}\leq \frac{\lambda_1(v)}{\Theta_1(v)}.\]

This completes the proof of Lemma \ref{3lam}.
\end{proof}

\begin{itemize}
	\item {\bf Criterion 2:}    The sequence $\{\Lambda_{II}(\kappa)\}_{\kappa\in \mathbb{N}}$ is  bounded.
\end{itemize}

Theorem \ref{mathe} will be proved by verifying this criterion.   As warm-up, we consider $H_{n_0}$ for a fixed $n_0$ large enough.
\begin{Lemma}\label{initi}
There exists a constant $N_0 > 0$ such that for any $n_0 > N_0$, if $v \in (\kappa_0)_{II}$ with $\kappa_0 = \phi(n_0) \in S$, then $\Theta_1(v) \sim \lambda_1(v)$. Namely, there exists a constant $C > 0$, independent of $n$, such that
\[
\frac{1}{C} < \frac{\lambda_1(v)}{\Theta_1(v)} \leq C.
\]
\end{Lemma}
\begin{proof}
Note that $\lambda_1(v)\sim\|q_{n_0}\alpha\|\sim q_{n_0}^{-1}$. Based on (\ref{Z11}) and (\ref{Z13}), we write $\bar{N}:=\bar{N}(\kappa_0)$. It follows from Proposition \ref{eveco} that one can find $k_0,k_1\in \{0,\ldots,q_{\bar{N}-1}\}$ such that
\[\Theta_1(F^{k_0}(v))\gtrsim \lambda_1(v),\quad  \Theta_1(F^{k_1}(v))\lesssim\lambda_1(v).\]
Without loss of generality, we assume $k_0\leq k_1$.
Let $v$ be a chord connecting $x_i$ and $x_j$. Then $\theta(v)=x_j-x_i$. It follows that
\[\Theta_1(F^{k_0}(v))=\|x_{j+k_0}-x_{i+k_0}\|,\quad \Theta_1(F^{k_1}(v))=\|x_{j+k_1}-x_{i+k_1}\|.\]
By Definition \ref{kvec3}, $|j-i|\sim q_{n_0}$.
By Lemma \ref{disne}, for each $i\in \mathbb{Z}$, we have
\[|y_{i+1}-y_i|\lesssim \frac{1}{q_{n_0}^{3+\varepsilon}}.\]
By the construction of $G_n$,
\[x_{j+k_0}-x_{i+k_0}=x_{j+k_1}-x_{i+k_1}-\sum_{l=0}^{k_1-k_0-1}(y_{j+l}-y_{i+l}).\]
Note that $0\leq k_1-k_0\leq q_{\bar{N}-1}$. Then
\begin{equation}\label{i-jrole}
\left|\sum_{l=0}^{k_1-k_0-1}(y_{j+l}-y_{i+l})\right|\lesssim \frac{1}{q_{n_0}^{1+\frac{\varepsilon}{2}}},
\end{equation}
which implies
\[\|x_{j+k_1}-x_{i+k_1}\|\sim\|x_{j+k_0}-x_{i+k_0}\|\sim \frac{1}{q_{n_0}}\sim\lambda_1(v).\]
Moreover, we obtain  \[\Theta_1(v)=\|x_j-x_i\|\sim\frac{1}{q_{n_0}}\sim \lambda_1(v).\]
\end{proof}

\subsection{Reduction}
Let us recall  $v$ connects $(x_i,y_i)$ and  $(x_j,y_j)$. Let $s(v)$ be the first-order difference quotient along the minimal configuration, i.e. $s(v):=\frac{r(v)}{\Theta_1(v)}$ where $r(v):=y_j-y_i$. Based on the Aubry-Mather theory, we know that there exists a constant $s_{\max}$ such that $|s(v)|\leq s_{\max}$ (\cite[Theorem 14.1]{FoM}).  Note that one can take $s_{\max}=2$ for $G_n$.  Moreover, if
\[|\theta(v)|,\quad |\theta(F(v))|,\quad |\theta(F^{-1}(v))|\leq \frac{1}{2},\]
then
\begin{equation}\label{boud}
\frac{1}{3}\leq \frac{1}{1+s_{\max}}\leq \frac{\Theta_1(F(v))}{\Theta_1(v)}\leq 1+s_{\max}\leq 3.
\end{equation}

We also note that $\frac{\Theta_1(F^j(v))}{\Theta_1(v)}$ behaves like a multiplicative cocycle. Inspired by the theory of  the circle diffeomorphism, we consider the expansion
\[K^0(j|v):=\ln \frac{\Theta_1(F^j(v))}{\Theta_1(v)}.\]

Denote for $\iota\in \{I, II\}$,
\[{K^0(v):=K^0(1|v)};\]
\[\bold{K}_N^0(\kappa):=\sup\{|K^0(q_N|v)|:\ v\in (\kappa)_{\iota}\};\]
\begin{equation}\label{bigK}
\widetilde{\mathbf{K}}_N^0(\kappa):=\sup\{|K^0(j|v)|:\ 0\leq j\leq q_N,\ \ v\in (\kappa)_{\iota}\}.
\end{equation}
Then
\[K^0(-j|v)=-K^0(j|F^{-j}(v)).\]
It follows that
\begin{equation}\label{skno}
\widetilde{\mathbf{K}}_N^0(\kappa)=\sup\{|K^0(j|v)|:\ -q_N\leq j\leq q_N,\ \ v\in (\kappa)_{\iota}\}.
\end{equation}

\begin{Lemma}\label{Cc} Assume $v\in (\kappa)_{II}$.
If there exists a constant $C$ independent of $\kappa$ such that for each $\kappa\in \mathbb{N}$,
\begin{equation}\label{ttkk}
\widetilde{\mathbf{K}}_{\bar{N}(\kappa)}^0(\kappa)\leq C,
\end{equation}
then {\bf Criterion 2} holds, where $\bar{N}(\kappa)$ is defined by (\ref{Z11}).
\end{Lemma}
\begin{proof} We use the same strategy as the proof of Lemma \ref{initi}. Note that for $\kappa$ sufficiently large such that $n_{\kappa} > N_0$, where $N_0$ is given by Lemma~\ref{initi}, we have $\lambda_1(v) \sim \|q_{n_\kappa} \alpha\| \sim q_{n_\kappa}^{-1}$.
By Proposition \ref{eveco}, one can find $k_0,k_1\in \{0,\ldots,q_{\bar{N}(\kappa)-1}\}$ such that
\[\Theta_1(F^{k_0}(v))\gtrsim \lambda_1(v),\quad  \Theta_1(F^{k_1}(v))\lesssim\lambda_1(v).\]
Given an index $j\in \{0,\ldots,q_{\bar{N}(\kappa)-1}\}$, one has either
\begin{equation}\label{eithor}
\Theta_1(F^{j}(v))\gtrsim \lambda_1(v)\ \ \mathrm{or}\ \ \Theta_1(F^{j}(v))\lesssim\lambda_1(v).
\end{equation}
We assume the first case happens.
Due to $|k_1-j|\leq q_{\bar{N}(\kappa)}$, it follows from (\ref{ttkk}) that
\[e^{-C}\leq \frac{\Theta_1(F^{k_1}(v))}{\Theta_1(F^{j}(v))}\leq e^C.\]
It yields
\[\Theta_1(F^{j}(v))\sim\Theta_1(F^{k_1}(v))\sim \lambda_1(v).\]
For the second case in (\ref{eithor}), we obtain the same result by comparing $\Theta_1(F^{j}(v))$ with $\Theta_1(F^{k_0}(v))$.
In particular, we have $\Theta_1(v)\sim \lambda_1(v)$.
This completes the proof.
\end{proof}

Therefore, in order to prove Criterion 2 holds, we only need to verify
\begin{itemize}
	\item {\bf Criterion 3:}   The sequence  $\{\widetilde{\mathbf{K}}_{\bar{N}(\kappa)}^0(\kappa)\}_{\kappa\in \mathbb{N}}$ is bounded.
\end{itemize}

\vspace{2em}

\section{Verification of Criterion 3 for a given $\kappa$}\label{se4}
Following \cite{KO1}, we introduce some notations  to describe the minimal configuration with frequency $\alpha$.
\subsection{Distortions}
Let $v_1,v_2$ be two  chords. Denote by $\lambda_2(v_1,v_2)$ the average of the $\lambda$-length of the span of the two chords and that of the part  between them. To verify Criterion 3, we only need to consider the case with $\lambda_2(v_1,v_2)$ small enough. We assume, once and for all, $v_1$ connects $(x_{i_1},y_{i_1})$ and  $(x_{j_1},y_{j_1})$ and  $v_2$ connects $(x_{i_2},y_{i_2})$ and  $(x_{j_2},y_{j_2})$. Moreover, we require
\[0\leq x_{i_1}<x_{j_1}\leq x_{i_2}<x_{j_2}\leq 1,\quad 0<x_{j_2}-x_{i_1}\leq \frac{1}{2}. \]
Then
\[\lambda_2(v_1,v_2)=\frac{1}{2}(\|(i_2-i_1)\alpha\|+\|(j_2-j_1)\alpha\|).\]
We also denote
\[\theta(v_1,v_2):=x_{i_2}-x_{i_1},\quad \Theta_2(v_1,v_2):=\|\theta(v_1,v_2)\|.\]
By  the orientation preserving property (see Proposition \ref{OP}), we have
\begin{equation}\label{thev1}
\Theta_2(v_1,v_2)\geq \Theta_1(v_1).
\end{equation}
If  $v_1,v_2$ are two  chords satisfying $v_i\in (\kappa)_{II}$ for $i=1,2$, then $|i_2-i_1|\sim |j_2-j_1|$.

We denote $(v_1,v_2)\in (\kappa,r)_{II}$ if   $v_i\in (\kappa)_{II}$ for $i=1,2$ and
\[\frac{1}{q_{r+2}}\leq\lambda_2(v_1,v_2)\leq \frac{1}{q_{r}}.\]
Moreover, Proposition \ref{eveco} works.
Let \[K^1(j|v_1,v_2):=K^0(j|v_2)-K^0(j|v_1),\]
\[{K^0(v_1,v_2)=K^1(1|v_1,v_2)}.\]
The distortion of the orbits with the length $q_N$ generated by $F$ is defined by
\[\bold{K}_N^1(\kappa,r):=\sup\{|K^1(q_N|v_1,v_2)|:\ (v_1,v_2)\in (\kappa,r)_{II}\},\]
\begin{equation}\label{distK}
\widetilde{\mathbf{K}}_N^1(\kappa,r):=\sup\{|K^1(j|v_1,v_2)|:\ 0\leq j\leq q_N,\ \ (v_1,v_2)\in (\kappa,r)_{II}\}.
\end{equation}
Similar to (\ref{skno}), it is easy to obtain
\[\widetilde{\mathbf{K}}_N^1(\kappa,r)=\sup\{|K^1(j|v_1,v_2)|:\ -q_N\leq j\leq q_N,\ \ (v_1,v_2)\in (\kappa,r)_{II}\}.\]

By using Denjoy's idea on the circle diffeomorphism \cite{De}, it was shown that (\cite[Lemma 2.4]{KO1})
\begin{Proposition}\label{ko11}
For each $\kappa\in \mathbb{N}$, and each $N\leq \tilde{N}(\kappa)$, there holds \[\bold{K}_N^0(\kappa)\leq 2 \widetilde{\mathbf{K}}_N^1(\kappa,N).\]
\end{Proposition}
\subsection{Difference quotients}
Recall that $v_1$ connects $(x_{i_1},y_{i_1})$ and  $(x_{j_1},y_{j_1})$ and  $v_2$ connects $(x_{i_2},y_{i_2})$ and  $(x_{j_2},y_{j_2})$.
We also assume $v_3$ connects $(x_{i_3},y_{i_3})$ and  $(x_{j_3},y_{j_3})$ and  $v_4$ connects $(x_{i_4},y_{i_4})$ and  $(x_{j_4},y_{j_4})$. Additionally,
assume they satisfy
\[0\leq x_{i_1}<x_{j_1}\leq x_{i_2}<x_{j_2}\leq x_{i_3}<x_{j_3}\leq x_{i_4}<x_{j_4}\leq 1,\quad 0<x_{j_4}-x_{i_1}\leq \frac{1}{2}. \]
Define
\[\lambda_4(v_1,v_2,v_3,v_4):=\frac{1}{2}(\lambda(v_1,v_3)+\lambda(v_2,v_4)).\]
Accordingly, we denote
\[\theta(v_1,\ldots,v_4):=\frac{x_{j_4}+x_{i_4}+x_{j_3}+x_{i_3}-x_{j_2}-x_{i_2}-x_{j_1}-x_{i_1}}{4},\]
\[ \Theta_4(v_1,\ldots,v_4):=\|\theta(v_1,\ldots,v_4)\|.\]
If $(v_1,v_2), (v_3,v_4)\in (\kappa,r)_{II}$, then $|i_3-i_1|\sim |i_4-i_2|$.
We denote $(v_1,v_2,v_3,v_4)\in (\kappa,r,s)_{II}$ for $s\leq r$ if  they are four  chords satisfying  $(v_1,v_2), (v_3,v_4)\in (\kappa,r)_{II}$ and
\[\frac{1}{q_{s+2}}\leq\lambda_4(v_1,v_2,v_3,v_4)\leq\frac{1}{q_{s}}.\]

To analyze the dynamics and geometry of the minimal configuration, we have to consider the second-order and third-order difference quotients, denoted by
\[\nabla^1(v_1,v_2):=\frac{s(v_2)-s(v_1)}{\Theta_2(v_1,v_2)},\] \[\nabla^2(v_1,\ldots,v_4):=\nabla^1(v_3,v_4)-\nabla^1(v_1,v_2).\]
Let
\begin{equation}\label{na11}
\nabla^1(\kappa,r):=\sup_{(v_1,v_2)\in (\kappa,r)_{II}}|\nabla^1(v_1,v_2)|.
\end{equation}
\begin{equation}\label{na22}
\nabla^2(\kappa,r,s):=\sup_{(v_1,\ldots,v_4)\in (\kappa,r,s)_{II}}|\nabla^2(v_1,\ldots,v_4)|.
\end{equation}
We connect $\nabla^1(v_1,v_2)$ and the dynamics of $F$ in the following way.
Denote
\begin{equation}\label{Ef}
E^1(j|v_1,v_2):=\Theta_2(v_1,v_2)^{-1}\left(\frac{\Theta_1(F^j(v_2))}{\Theta_1(v_2)}-\frac{\Theta_1(F^j(v_2))}{\Theta_1(v_2)}\right).
\end{equation}
\begin{Lemma}\label{v1v2}
Let $(v_1,v_2)\in (\kappa,r)_{II}$ with $\kappa\geq n_0$, where $n_0$ satisfies Lemma \ref{initi}. Then for each $j\in\mathbb{Z}$,
\[E^1(1|F^j(v_1),F^j(v_2))= \nabla^1(F^j(v_1),F^j(v_2)).\]
\end{Lemma}
\begin{proof} Note that $\lambda$ is an invariant metric with respect to $F$.   For each $j\in\mathbb{Z}$, we have
\[\lambda_1(F^j(v_i))=\lambda(v_i)\  (i=1,2), \quad \lambda_2(F^j(v_1),F^j(v_2))=\lambda_2(v_1,v_2).\]It follows that $(F^j(v_1),F^j(v_2))\in (\kappa,r)$ for each $j\in\mathbb{Z}$. Without loss of generality, we fix $j=0$. It suffices to prove
\begin{equation}\label{lef}
E^1(1|v_1,v_2)= \nabla^1(v_1,v_2).
\end{equation}
For a given $v\in (\kappa)_{II}$,  we assume $\theta(v)=x_{l+q_\kappa}-x_l$.
Since $\alpha$ is irrational, then by the orientation-preserving property (see Proposition \ref{OP}) and Lemma \ref{initi},
\[\Theta_1(v)=\|x_{l+q_\kappa}-x_{l}\|<\|x_{l+q_{n_0}}-x_{l}\|\sim \frac{1}{q_{n_0}},\quad \mathrm{for}\ \kappa\geq n_0.\]

By the definition of $H_{n_0}$,
\begin{align*}
	&\frac{\Theta_1(F(v_2))}{\Theta_1(v_2)}-\frac{\Theta_1(F(v_1))}{\Theta_1(v_1)}\\
	=&\frac{x_{l_2+q_\kappa}-x_{l_2}+y_{l_2+q_\kappa}-y_{l_2}}{x_{l_2+q_\kappa}-x_{l_2}}-
	\frac{x_{l_1+q_\kappa}-x_{l_1}+y_{l_1+q_\kappa}-y_{l_1}}{x_{l_1+q_\kappa}-x_{l_1}}\\
	=&s(v_2)-s(v_1),
\end{align*}
which implies the left hand side of (\ref{lef}) is nothing but
\[\frac{s(v_2)-s(v_1)}{\Theta_2(v_1,v_2)}=\nabla^1(v_1,v_2).\]
This completes the proof of Lemma \ref{v1v2}.
\end{proof}
\subsection{Growth of the distortions}
We assert the following
\begin{Lemma}\label{last}
Given $\kappa\in \mathbb{N}$, assume
\begin{equation}\label{f1x}
\nabla^1(\kappa,\tilde{N})\leq \frac{\varepsilon}{960A_\alpha},
\end{equation}
\begin{equation}\label{f2x}
\nabla^2(\kappa,\tilde{N},\tilde{N})\leq C_0 q_{\tilde{N}}^{-\frac{\varepsilon}{3}},
\end{equation}
where $\tilde{N}:=\tilde{N}(\kappa)$ is defined by (\ref{Z12}), $\varepsilon$ is given by (\ref{epss}) and $C_0$ is a constant independent of $\kappa$. Then for each $0\leq r\leq \tilde{N}$, we have
\begin{equation}\label{kkes}
\widetilde{\mathbf{K}}_{r}^1(\kappa,r)\lesssim (\frac{9}{10})^{\frac{\varepsilon}{2} r}.
\end{equation}
\end{Lemma}
For consistency, we first verify Criterion 3 by using the conclusion of Lemma \ref{last}.
Let us recall $\bar{N}:=\bar{N}(\kappa)$ is defined by (\ref{Z11}). For a given $j$ with $0\leq j<q_{\bar{N}}$, it can be written as
\[j=b_Rq_R+\ldots+b_rq_r+\ldots+b_0q_0,\]
where $q_0=1$, $q_R\leq j<q_{R+1}$ and $0\leq b_r\leq a_{r+1}$. In particular,
\[j=c q_{\tilde{N}}+\sum_{r=0}^{\tilde{N}-1}b_rq_r,\]
where $c\leq q_{\bar{N}}q_{\tilde{N}}^{-1}$. A direct calculation shows that for $v\in (\kappa)_{II}$,
\begin{equation}\label{kese}
K^0(j|v)\leq a_{R+1}\bold{K}_R^0(\kappa)+\ldots+a_{r+1}\bold{K}_{r}^0(\kappa)+\ldots+a_1\bold{K}_0^0(\kappa).
\end{equation}
It follows from (\ref{kese}) and Proposition \ref{ko11} that
\begin{equation}\label{nrev}
	\begin{split}
		\widetilde{\mathbf{K}}^0_{\bar{N}}(\kappa)&\leq c\bold{K}^0_{\tilde{N}}(\kappa)+\sum_{r=0}^{\tilde{N}-1}a_r\bold{K}^0_{r}(\kappa)\\
		&\leq 2c\widetilde{\mathbf{K}}^1_{\tilde{N}}(\kappa,\tilde{N})+2\sum_{r=0}^{\tilde{N}-1}a_r\widetilde{\mathbf{K}}^1_{r}(\kappa,r).
	\end{split}
\end{equation}
Note that for a given $\kappa\in \mathbb{N}$, by (\ref{Z11}) and (\ref{Z12}),
\begin{equation}\label{qratq}
c\leq \frac{q_{\bar{N}}}{q_{\tilde{N}}}=\frac{q_{n_\kappa+2\gamma_0}}{q_{n_\kappa-2\gamma_0}}\sim 1.
\end{equation}

\noindent Therefore, we obtain
\begin{equation}\label{nrevxx}
	\begin{split}
		\widetilde{\mathbf{K}}^0_{\bar{N}}(\kappa)&\lesssim \widetilde{\mathbf{K}}^1_{\tilde{N}}(\kappa,\tilde{N})+\sum_{r=0}^{\tilde{N}-1}\widetilde{\mathbf{K}}^1_{r}(\kappa,r),\\
		&\lesssim  (\frac{9}{10})^{\frac{\varepsilon}{2} \tilde{N}}+\ldots+(\frac{9}{10})^{\frac{\varepsilon}{2} r}+\ldots+1 ,\\
		&\lesssim \sum_{r=0}^{+\infty}(\frac{9}{10})^{\frac{\varepsilon}{2} r},
	\end{split}
\end{equation}
which implies {$\widetilde{\mathbf{K}}_{\bar{N}(\kappa)}^0(\kappa)\lesssim 1$.}

\subsection{Proof of Lemma \ref{last}}

To help the reader follow the argument, we first outline the main steps.

\begin{itemize}
	\item[(1)] Lemma~\ref{A1} provides monotonicity properties of \(\nabla^1(\kappa,r)\) with respect to \(\kappa\) and \(r\).
	
	\item[(2)] The proof proceeds in two stages. Lemma~\ref{keyy1} supplies three basic estimates that are needed to establish (\ref{kkes}); none of them uses condition (\ref{f2x}). Specifically,
	\begin{itemize}
		\item Lemma~\ref{keyy1}(1) gives a rough upper bound for \(\widetilde{\mathbf{K}}_{r}^1(\kappa,r)\);
		\item Lemma~\ref{keyy1}(2) is used to obtain (\ref{lem6i}) in the proof of Lemma~\ref{karrr1}; for a more compact presentation we include it as part of Lemma~\ref{keyy1}(2);
		\item Lemma~\ref{keyy1}(3) follows from a direct computation and clarifies under which relation between \(r\) and \(R\) (roughly \(r\ll R\)) the rough bound of Lemma~\ref{keyy1}(1) actually yields exponential decay of \(\widetilde{\mathbf{K}}_{r}^1(\kappa,R)\).
	\end{itemize}
	
	\item[(3)] Lemma~\ref{keyy1} shows that when \(r\ll R\), one already obtains exponential decay of \(\widetilde{\mathbf{K}}_{r}^1(\kappa,R)\). Lemma~\ref{karrr1} proves that, with the aid of condition (\ref{f2x}), this exponential decay persists even when \(r=R\).
\end{itemize}

\subsubsection{Monotonicity of $\nabla^1(\kappa,r)$}

Let $v\in (\kappa)_{\iota}$ with $\iota\in \{I, II\}$ for $\kappa$ large enough. Note that if $v=v'+v''$, then $\Theta_1(F^j(v))=\Theta_1(F^j(v'))+\Theta_1(F^j(v''))$ for each $j\in \mathbb{Z}$. Moreover,
\[\frac{\Theta_1(F^{j}(v))}{\Theta_1(v)}=\frac{\Theta_1(v')}{\Theta_1(v)}\frac{\Theta_1(F^{j}(v'))}{\Theta_1(v')}+\frac{\Theta_1(v'')}{\Theta_1(v)}
\frac{\Theta_1(F^{j}(v''))}{\Theta_1(v'')}.\]
Due to the convexity of  $\ln x$ with respect to $x$, it follows that for any $\kappa\in\mathbb{N}$,
\begin{equation}\label{kinc}
\bold{K}_r^0(\kappa)\leq \bold{K}_r^0(\kappa+1).
\end{equation}

\begin{Lemma}\label{A1}
For any $r\leq \kappa$, we have
\[\nabla^1(\kappa,r)\leq \nabla^1(\kappa+1,r)\leq \nabla^1(\kappa+1,r+1).\]
\end{Lemma}

\begin{proof} We recall that \( s(v)=\frac{r(v)}{\Theta_1(v)} \). For any two chords \( v \) and \( v' \), suppose \( v \) can be decomposed as \( v = v_1 + v_2 \). From the convex combination
\[
s(v_1+v_2) = \frac{\Theta_1(v_1)}{\Theta_1(v_1)+\Theta_1(v_2)}\,s(v_1) + \frac{\Theta_1(v_2)}{\Theta_1(v_1)+\Theta_1(v_2)}\,s(v_2),
\]
and by direct computation we obtain the decomposition of \( \nabla^1(v_1+v_2,v') \):
\begin{align}
	\nabla^1(v,v') &= \frac{s(v') - b_1 s(v_1) - b_2 s(v_2)}{\Theta(v_1+v_2,v')} \notag \\
	&= \frac{\Theta_2(v_1,v')}{\Theta_2(v_1+v_2,v')}\,b_1\nabla^1(v_1,v') + \frac{\Theta_2(v_2,v')}{\Theta_2(v_1+v_2,v')}\,b_2\nabla^1(v_2,v') \notag \\
	&= c_1\nabla^1(v_1,v') + c_2\nabla^1(v_2,v'),
\end{align}
where the coefficients are given by
\[
b_1 = \frac{\Theta_1(v_1)}{\Theta_1(v_1)+\Theta_1(v_2)},\qquad
b_2 = \frac{\Theta_1(v_2)}{\Theta_1(v_1)+\Theta_1(v_2)},\]
\[
c_1 = \frac{\Theta_2(v_1,v')}{\Theta_2(v_1+v_2,v')}b_1,\qquad
c_2 = \frac{\Theta_2(v_2,v')}{\Theta_2(v_1+v_2,v')}b_2.
\]
A direct calculation shows that \( c_1+c_2 = b_1+b_2 = 1 \).

More generally, suppose \( v \) and \( v' \) can be split as
\[
v = \sum_{i=1}^{I} v_i,\qquad v' = \sum_{j=1}^{J} v_j',
\]
with \( I,J\in\mathbb{N} \). Then
\[
\nabla^1(v,v') = \sum_{i=1}^{I} d_i\Bigl(\sum_{j=1}^{J} e_j\nabla^1(v_i,v_j')\Bigr),
\]
i.e., \( \nabla^1(v,v') \) can be expressed as a convex combination of the quantities \( \nabla^1(v_i,v_j') \).

{Using the arithmetic properties of a rotation number of constant type, for given \( v,v'\in (\kappa,r)_{II} \) we can split \( v \) and \( v' \) into several chords such that each \( v_i, v_j'\in (\kappa+1,r)_{II} \).} By definition we then obtain
\[
\nabla^1(\kappa,r) \le \nabla^1(\kappa+1,r).
\]

By Lemma~\ref{v1v2}, if \( (v_1,v_2)\in (\kappa+1,r)_{II} \), then
\[
\nabla^1(v_1,v_2)=\Theta_2(v_1,v_2)^{-1}\Bigl(\frac{\Theta_1(F(v_2))}{\Theta_1(v_2)}-\frac{\Theta_1(F(v_2))}{\Theta_1(v_2)}\Bigr).
\]
{Again using the arithmetic properties of a constant-type rotation number, we can insert a chord \( v_3\in (\kappa+1)_{II} \) between \( v_1 \) and \( v_2 \) so that
	\[
	(v_1,v_3)\in (\kappa+1,r+1)_{II},\qquad (v_3,v_2)\in (\kappa+1,r+1)_{II}.
	\]}
Note that
\[
\nabla^1(v_1,v_2)=\frac{\Theta_2(v_1,v_3)}{\Theta_2(v_1,v_2)}\,\nabla^1(v_1,v_3)+\frac{\Theta_2(v_3,v_2)}{\Theta_2(v_1,v_2)}\,\nabla^1(v_3,v_2),\]
\[\Theta_2(v_1,v_2)=\Theta_2(v_1,v_3)+\Theta_2(v_3,v_2).
\]
Hence, by definition,
\[
\nabla^1(\kappa+1,r)\le \nabla^1(\kappa+1,r+1).
\]
\end{proof}

\subsubsection{Exponential decay of $\widetilde{\mathbf{K}}_r^1(\kappa,R)$ with $r\ll R$}
For simplicity of notation, we write
\[
F^i(v_1,\dots,v_4) \;:=\; \bigl(F^i(v_1),\dots,F^i(v_4)\bigr),
\]
and we adopt analogous conventions for the action of \(F\) on pairs of chords.

Denote
\[K^2(v_1,\ldots,v_4):=\frac{K^1(v_3,v_4)}{\Theta_2(v_3,v_4)}-\frac{K^1(v_1,v_2)}{\Theta_2(v_1,v_2)};\]
\[K^2(j|v_1,\ldots,v_4):=\frac{K^1(j|v_3,v_4)}{\Theta_2(v_3,v_4)}-\frac{K^1(j|v_1,v_2)}{\Theta_2(v_1,v_2)}.\]
It follows that
\begin{equation}\label{sk2}
\Theta_2(v_3,v_4)K^2(j|v_1,\ldots,v_4)=K^1(j|v_3,v_4)-\frac{\Theta_2(v_3,v_4)}{\Theta_2(v_1,v_2)}K^1(j|v_1,v_2).
\end{equation}
A direct calculation shows (see \cite[Proposition 2.3]{KO1})
\begin{Lemma}\label{K22}
Given $(v_1,v_2,v_3,v_4)\in (\kappa,r,s)_{II}$, we have
\[K^2(v_1,\ldots,v_4)\lesssim \nabla^2(v_1,\ldots,v_4)+\Theta_4(v_1,\ldots,v_4).\]
\end{Lemma}

Before proving the exponential decay of $\widetilde{\mathbf{K}}_R^1(\kappa,R)$, we first present a rough estimate that does not rely on condition (\ref{f2x}). Building on this, we then use (\ref{f2x}) to obtain the finer estimate in Lemma~\ref{last}.
\begin{Lemma}\label{keyy1}
Given $\kappa\in\mathbb{N}$, and $0\leq r\leq R\leq \tilde{N}$ where  $\tilde{N}:=\tilde{N}(\kappa)$ is defined by (\ref{Z12}).
\begin{itemize}
	\item [(1)] There holds $\widetilde{\mathbf{K}}_r^1(\kappa,R)\lesssim \varepsilon e^{\frac{\varepsilon}{40}R}q_rq_R^{-1}$.
	\item [(2)] For $j\leq q_r$ and $(v_1,\ldots,v_4)\in (\kappa,R,R)_{II}$, we have
	\begin{align*}
		K^2(j|v_1,\ldots,v_4)-\sum_{i=0}^{j-1}K^2(F^i(v_1,\ldots,v_4))
		\lesssim \varepsilon e^{\frac{\varepsilon}{10}R}q_rq_R^{-1}\sum_{i=0}^{j-1}\frac{\Theta_2(F^i(v_3,v_4))}{\Theta_2(v_3,v_4)}.
	\end{align*}
	\item [(3)] In particular, if $q_R^{1-\frac{\varepsilon}{2}}\sim q_r$, then
	\begin{equation}\label{eeps}
	e^{\frac{\varepsilon}{40}R}q_rq_R^{-1}\leq e^{\frac{\varepsilon}{10}R}q_rq_R^{-1}\sim e^{\frac{\varepsilon}{10}R}q_R^{-\frac{\varepsilon}{2}}<e^{\frac{\varepsilon}{10}R}q_R^{-\frac{\varepsilon}{3}}\leq (\frac{9}{10})^{\frac{\varepsilon R}{2}}.
	\end{equation}
\end{itemize}
\end{Lemma}
\begin{proof}
By (\ref{D2}), $q_R\geq \sqrt{2}^{R}$. For Item (3), a direct calculation shows
\[e^{\frac{\varepsilon}{10}R}q_R^{-\frac{\varepsilon}{2}}\leq e^{\frac{\varepsilon}{10}R}(\frac{1}{\sqrt{2}})^{\frac{\varepsilon}{2}R}\leq (\frac{9}{10})^{\frac{\varepsilon R}{2}}.\]

For $(v_1,v_2)\in (\kappa,r)_{II}$ with $0\leq r\leq R$, denote by $J(v_1,v_2)$ the arc on $\mathbb{T}$ spanned by $v_1\cup v_2$. By definition, we have
\[\frac{1}{q_{r+2}}<\lambda_2(v_1,v_2)<\frac{1}{q_r}.\]
Then for $0\leq j<q_r$, no point on $\mathbb{T}$ is  contained in more than two of the arcs $\{F^i(J(v_1,v_2))\}_{i=0}^j$ (see \cite[p65, Proposition (8.3)]{H11}). It implies
\begin{equation}\label{hcir}
\sum_{i=0}^j\Theta_2(F^i(v_1,v_2))\leq 2.
\end{equation}
By definition, we have
\begin{align*}
	K^1(j|v_1,v_2)&=K^0(j|v_2)-K^0(j|v_1)\\
	&=\sum_{i=0}^{j-1}\left(\ln \frac{\Theta_1(F^{i+1}(v_2))}{\Theta_1(F^i(v_2))}-\ln \frac{\Theta_1(F^{i+1}(v_1))}{\Theta_1(F^i(v_1))}\right).
\end{align*}
By (\ref{boud}),  $\frac{\Theta_1(F^{i+1}(v))}{\Theta_1(F^i(v))}$ is bounded for each $v\in (\kappa)_{II}$ and $i\in \mathbb{Z}$. It follows from the mean value theorem and Lemma \ref{v1v2} that for $0\leq j<q_r$,
\begin{equation}\label{KIJ}
|K^1(j|v_1,v_2)|\leq3\nabla^1(\kappa,r)\sum_{i=0}^j\Theta_2(F^i(v_1,v_2))\leq \frac{\varepsilon}{80A_\alpha}.
\end{equation}
Note that \begin{equation}\label{Rgam}
r\leq R+2\gamma_0\leq \tilde{N}(\kappa)+2\gamma_0=n_{\kappa}.
\end{equation}
It follows from   (\ref{f1x}) and Lemma \ref{A1} that
\[\widetilde{\mathbf{K}}_{r}^1(\phi(r),r)\leq \frac{\varepsilon}{80A_\alpha},\quad\widetilde{\mathbf{K}}_{r}^1(\kappa,r)\leq \frac{\varepsilon}{80A_\alpha},\]
where $\phi$ is defined by (\ref{phiex}).
By Proposition \ref{ko11} and (\ref{kese}), for $0<j<q_{R+2\gamma_0}$, where $\gamma_0$ is defined by (\ref{Z11}), we have
\[K^0(j|v)\leq 2a_{R+2\gamma_0}\widetilde{\mathbf{K}}_{R+2\gamma_0-1}^1(\kappa,R+2\gamma_0-1)+\ldots+2a_{l+1}\widetilde{\mathbf{K}}_{l}^1(\kappa,l)+\ldots+2a_1\widetilde{\mathbf{K}}_0^1(\kappa,0),\]
where $1\leq a_l\leq A_\alpha$. Similar to the argument in (\ref{nrev}), (\ref{qratq}) and (\ref{nrevxx}), we have
\begin{equation}\label{newb}
\widetilde{\mathbf{K}}_{R+2\gamma_0}^0(\kappa)\leq \frac{\varepsilon (R+c)}{40},
\end{equation}
where $0<c\leq q^{-1}_{R}q_{R+2\gamma_0}$. Due to the arithmetics of $\alpha$, $c$ is a constant independent of $R$.

By definition, for $(v_1,v_2)\in (\kappa,R)_{II}$, there exists a chord $\bar{v}$ determined by $v_1,v_2$ such that $\lambda_1(\bar{v})=\lambda_2(v_1,v_2)$ and $\bar{v}\in ({\phi(R)})_{I}$.  Note from (\ref{kinc}) that
\[\widetilde{\mathbf{K}}_{R+2\gamma_0}^0({\phi(R)})\leq\widetilde{\mathbf{K}}_{R+2\gamma_0}^0(\kappa)\leq \frac{\varepsilon (R+c)}{40}.\]
Since $\bar{v}\in ({\phi(R)})_{I}$, then $\lambda_1(\bar{v})\geq \frac{2}{q_{R+2\gamma_0}}$. By Proposition \ref{eveco}, there exists $0\leq j_0<q_{R+2\gamma_0}$ such that
\[\Theta_1(F^{j_0}(\bar{v}))\lesssim \lambda_1(\bar{v}).\]
By the definition of $\widetilde{\mathbf{K}}_{R+2\gamma_0}^0({\phi(R)})$ (see (\ref{bigK}) and (\ref{skno})), for each $j$ with $|j-j_0|<q_{R+2\gamma_0}$, we have
\[-\frac{\varepsilon (R+c)}{40}\leq \ln\frac{\Theta_1(F^j(\bar{v}))}{\Theta_1(F^{j_0}(\bar{v}))}\leq \frac{\varepsilon (R+c)}{40}.\]
Note that for each $i\in\mathbb{Z}$,
\[F^i(\bar{v})=F^i(v_1,v_2),\quad \lambda_1(\bar{v})=\lambda_2(F^i(v_1,v_2)).\]
It follows that for each $j$ with $|j-j_0|<q_{R+2\gamma_0}$, we have
\[\frac{\Theta_2(F^j(v_1,v_2))}{\lambda_2(F^{j}(v_1,v_2))}=\frac{\Theta_1(F^j(\bar{v}))}{\lambda_1(\bar{v})}\lesssim\frac{\Theta_1(F^j(\bar{v}))}
{\Theta_1(F^{j_0}(\bar{v}))}\leq e^{ \frac{\varepsilon (R+c)}{40}}.\]
Then for $j<q_r$, it follows from (\ref{KIJ}) that
\[K^1(j|v_1,v_2)\lesssim \varepsilon  \sum_{i=0}^{q_r}\Theta_2(F^i(v_1,v_2))\leq \varepsilon e^{ \frac{\varepsilon (R+c)}{40}}\sum_{i=0}^{q_r}\lambda_2(F^i(v_1,v_2))\lesssim \varepsilon e^{\frac{\varepsilon}{40}R}q_rq_R^{-1}.\]
Hence, we have
\begin{equation}\label{B11}
\widetilde{\mathbf{K}}_{r}^1(\phi(R),R)\lesssim \varepsilon e^{\frac{\varepsilon}{40}R}q_rq_R^{-1},\quad\widetilde{\mathbf{K}}_{r}^1(\kappa,R)\lesssim \varepsilon e^{\frac{\varepsilon}{40}R}q_rq_R^{-1}.
\end{equation}
This completes the proof of Item (1).

To prove Item (2), we write from (\ref{sk2}) that
\begin{align*}
	\Theta_2(v_3,v_4)K^2(j|v_1,\ldots,v_4)=&\sum_{i=0}^{j-1}\left(K^1(F^i(v_3,v_4))-\frac{\Theta_2(F^i(v_3,v_4))}{\Theta_2(F^i(v_1,v_2))}K^1(F^i(v_1,v_2))\right)\\
	&-\sum_{i=0}^{j-1}\left(\frac{\Theta_2(v_3,v_4)}{\Theta_2(v_1,v_2)}-\frac{\Theta_2(F^i(v_3,v_4))}{\Theta_2(F^i(v_1,v_2))}\right)K^1(F^i(v_1,v_2)).
\end{align*}
It is clear to see
\begin{align*}
	&\frac{\Theta_2(F^i(v_1,v_2))}{\Theta_2(F^i(v_3,v_4))}\left(\frac{\Theta_2(v_3,v_4)}{\Theta_2(v_1,v_2)}-
	\frac{\Theta_2(F^i(v_3,v_4))}{\Theta_2(F^i(v_1,v_2))}\right)\\
	\leq &\left|\frac{\Theta_2(v_3,v_4)}{\Theta_2(F^i(v_3,v_4))}-\frac{\Theta_2(v_1,v_2)}{\Theta_2(F^i(v_1,v_2))}\right|
	\frac{\Theta_2(F^i(v_1,v_2))}{\Theta_2(v_1,v_2)}.
\end{align*}
By the definition of $\Theta_2(v_i,v_j)$, we have
\[\frac{\Theta_2(F^i(v_1,v_2))}{\Theta_2(v_1,v_2)}\leq e^{\frac{\varepsilon}{40}R},\quad\left|\frac{\Theta_2(v_3,v_4)}{\Theta_2(F^i(v_3,v_4))}
-\frac{\Theta_2(v_1,v_2)}{\Theta_2(F^i(v_1,v_2))}\right|\lesssim e^{\frac{\varepsilon}{40}R}\widetilde{\mathbf{K}}_r^1(\phi(R),R).\]
We have
\[\frac{\Theta_2(F^i(v_1,v_2))}{\Theta_2(F^i(v_3,v_4))}\left(\frac{\Theta_2(v_3,v_4)}{\Theta_2(v_1,v_2)}-
\frac{\Theta_2(F^i(v_3,v_4))}{\Theta_2(F^i(v_1,v_2))}\right)\lesssim \varepsilon e^{\frac{\varepsilon}{10}R}q_rq_R^{-1}\leq\varepsilon (\frac{9}{10})^{\frac{\varepsilon R}{2}}.\]
By Lemma \ref{v1v2} and the assumption on $\nabla^1(\kappa,r)$, we know that $\frac{K^1(F^i(v_1,v_2))}{\Theta_2(F^i(v_1,v_2))}$ is bounded. Combining with (\ref{sk2}), we obtain Item (2).
\end{proof}

\subsubsection{Exponential decay of $\widetilde{\mathbf{K}}_r^1(\kappa,R)$ with $r=R$}

Given $\kappa\in \mathbb{N}$ and $r\leq R\leq\tilde{N}(\kappa)$, we consider the quantity
\[\Delta(r;R,\kappa):={\sup}\left\{\frac{|K^1(q_r|v_1,v_2)|}{\sum_{i=0}^{q_r-1}\Theta_2(F^i(v_1,v_2))}\ \big|\ (v_1,v_2)\in (\kappa,R)_{II}\right\}.\]
To improve the estimate (\ref{B11}), the key point is to refine the estimate for (\ref{KIJ}) to an exponential one. Essentially, this improvement is guaranteed by condition (\ref{f2x}). Observe that in (\ref{KIJ}) one has (see (\ref{hcir}))
\begin{equation}\label{r222}
\sum_{i=0}^{q_r-1}\Theta_2\bigl(F^i(v_1,v_2)\bigr)\leq 2.
\end{equation}
Thus it suffices to prove that \(\Delta(r;R,\kappa)\) decays exponentially.

We provide another proposition on the arithmetics of $\alpha$ here, which is useful to prove Lemma \ref{karrr1} below.
\begin{Proposition}\label{arthh}
Given $\delta\in (0,\frac{\ln 2}{5\ln(1+A_\alpha)})$, for any $m,M\in \mathbb{N}$ with $M\geq m\geq 20$, if $q_M\leq q_m^{1+\delta}$, then $M\leq m+[\frac{m}{2}]\leq \frac{3}{2}m$.
\end{Proposition}
\begin{proof}
It suffices to prove $q_M\leq q_{m+[\frac{m}{2}]}$.  Due to the requirement on $\delta$, a direct calculation shows for $m\geq 20$,
\[(1+A_\alpha)^{\delta m}\leq 2^{\frac{m}{5}}\leq 2^{[\frac{m}{4}]}.\]
By using the arithmetics of $\alpha$, we have $q_m\leq (1+A_\alpha)^m$. Then by (\ref{D2}), we have
\[q_{m+[\frac{m}{2}]}\geq 2q_{m+[\frac{m}{2}]-2}\geq \ldots \geq 2^{[\frac{m}{4}]}q_m\geq q_m^{1+\delta},\]
which implies $q_M\leq q_{m+[\frac{m}{2}]}$.
\end{proof}

Based on Lemma \ref{keyy1}(1),  under the assumption (\ref{f1x}), there holds $\widetilde{\mathbf{K}}_r^1(\kappa,R)\lesssim \varepsilon e^{\frac{\varepsilon}{40}R}q_rq_R^{-1}$. It follows that for $q_R^{1-\frac{\varepsilon}{2}}\sim q_r$,
\begin{equation}\label{KES}
\widetilde{\mathbf{K}}_r^1(\kappa,R)\lesssim (\frac{9}{10})^{\frac{\varepsilon R}{2}},
\end{equation}
where $r\ll R$ due to the requirement $q_R^{1-\frac{\varepsilon}{2}}\sim q_r$. In the following, we will show that with the help of the assumption (\ref{f2x}), $\widetilde{\mathbf{K}}_r^1(\kappa,R)$ is still exponentially small even for $r=R$, i.e.
\[\widetilde{\mathbf{K}}_R^1(\kappa,R)\lesssim (\frac{9}{10})^{\frac{\varepsilon R}{2}}.\]
\begin{Lemma}\label{karrr1} If $q_R^{1-\frac{\varepsilon}{2}}\sim q_r$, then for $r$ large enough,
\begin{equation}\label{karrr}
\Delta(r;R,\kappa)\lesssim (\frac{9}{10})^{\frac{\varepsilon r}{2}}.
\end{equation}
Moreover,
\[\Delta(R;R,\kappa)\lesssim (\frac{9}{10})^{\frac{\varepsilon R}{2}},\quad \widetilde{\mathbf{K}}_R^1(\kappa,R)\lesssim (\frac{9}{10})^{\frac{\varepsilon R}{2}}.\]
\end{Lemma}
\begin{proof}
We prove this result by contradiction. We assume (\ref{karrr}) does not hold.  For fixed $R$ and $\kappa$, we write $\Delta(r)$ instead of $\Delta(r;R,\kappa)$ for simplicity.

At first, we claim that

\noindent{\bf Claim 1: }
For $2< r<R$, and $(v_1,v_2)\in (\kappa,R)_{II}$,
\begin{equation}\label{cla1}
\widetilde{\mathbf{K}}_R^1(\kappa,R)\leq 3\Delta(r).
\end{equation}
\noindent{\bf Proof of Claim 1: }
For $r<R$, $0\leq j<q_R$, break the run of length $j$ into blocks of length $q_r$ and a single remainder block. We obtain
\begin{equation}\label{wan}
{K^1(j|v_1,v_2)\leq \Delta(r)\sum_{i=0}^{j-1}\Theta_2(F^i(v_1,v_2))+\widetilde{\mathbf{K}}_r^1(\kappa,R).}
\end{equation}
If (\ref{karrr}) does not hold, it follows from (\ref{KES}) that the growth of \(\widetilde{\mathbf{K}}_r^1(\kappa,R)\) is negligible compared with \(\Delta(r)\). Together with (\ref{r222}) this yields
\[K^1(j|v_1,v_2)\leq 3\Delta(r),\]
which implies (\ref{cla1}).

\noindent{\bf Claim 2: }
\begin{equation}\label{kaka}
\Delta(R)\lesssim \Delta(r)^2+(\frac{9}{10})^{\frac{\varepsilon R}{2}}.
\end{equation}
\noindent{\bf Proof of Claim 2: }
Let $(v_1,v_2),(v'_1,v'_2)\in (\kappa,R)_{II}$ such that $K^1(q_R|v_1,v_2)$ and $K^1(q_R|v'_1,v'_2)$ have opposite signs. Without loss of generality, we assume \[K^1(q_R|v'_1,v'_2)\leq 0,\]otherwise, we exchange the roles of $v'_1$ and $v'_2$. Then for each $w'>0$, we have
\begin{equation}
|K^1(q_R|v_1,v_2)|\leq K^1(q_R|v_1,v_2)-w'K^1(q_R|v'_1,v'_2).
\end{equation}
Let $i_0\in (0,q_R-1)$ be such that $(F^{i_0}(v_1,v_2),v'_1,v'_2)\in (\kappa,R,R)_{II}$. It follows that for all $i\in \mathbb{Z}$,
\[(F^{i_0+i}(v_1,v_2),F^{i}(v'_1,v'_2))\in (\kappa,R,R)_{II}.\]
Denote
\[i':=\left\{\begin{array}{ll}
	\hspace{-0.4em}i-i_0,\quad i_0\leq i<q_R,\\
	\hspace{-0.4em} i-i_0+q_R,\quad 0\leq i<i_0.
\end{array}\right.\]
Then $F^{i_0+i}(v_1,v_2)$ and $F^{i}(v'_1,v'_2)$ are disjoint. We denote
\[Y(j|v_1,v_2):=\frac{\Theta_2(F^j(v_1,v_2))}{\Theta_2(v_1,v_2)},\]
{\[\Upsilon:=\frac{\Theta_2(F^{i'}(v'_1,v'_2))}{\Theta_2(v'_1,v'_2)}\cdot\frac{\Theta_2(F^{i_0}(v_1,v_2))}{\Theta_2(F^{i}(v_1,v_2))}\]}Then we have
\begin{equation}\label{longf}
\Upsilon=\left\{\begin{array}{ll}
	\hspace{-0.4em}\frac{Y(i'|v'_1,v'_2)}{Y(i'|F^{i_0}(v_1,v_2))},\quad i_0\leq i<q_R,\\
	\hspace{-0.4em} \frac{Y(i'|v'_1,v'_2)}{Y(i'|F^{i_0-q_R}(v_1,v_2))\cdot Y(q_R|F^{-q_R}(v_1,v_2))},\quad 0\leq i<i_0.
\end{array}\right.
\end{equation}
By taking the logarithm of the left hand side of (\ref{longf}) and combining with Proposition \ref{ko11}, we obtain
\[\Upsilon\in \left(\exp(-3\widetilde{\mathbf{K}}_R^1(\kappa,R)),\exp(3\widetilde{\mathbf{K}}_R^1(\kappa,R))\right).\]
By Claim 1, we know that
\[\Upsilon\in \left(\exp(-9\Delta(r)),\exp(9\Delta(r))\right).\]
Choosing
\[w'=\exp(-9\Delta(r))\frac{\Theta_2(F^{i_0}(v_1,v_2))}{\Theta_2(v'_1,v'_2)},\]
Then
\begin{equation}\label{midex}
\exp(-18\Delta(r))\leq w'\frac{\Theta_2(F^{i'}(v'_1,v'_2))}{\Theta_2(F^i(v_1,v_2))}\leq 1.
\end{equation}
Break the runs of length $i_0$ and $q_R-i_0$ into blocks of length $q_r$ plus one shorter block. Denote by $b_p$ the beginning of the $p'$th block. Write
\begin{equation}
\begin{split}
	K^1(q_R|v_1,v_2)-w'K^1(q_R|v'_1,v'_2)=\Xi_1+\Xi_2,
\end{split}
\end{equation}
where
\[ \Xi_1:=\sum_p\left(1-w'\frac{\Theta_2(F^{b'_p}(v'_1,v'_2))}{\Theta_2(F^{b_p}(v_1,v_2))}\right)K^1(q_r|F^{b_p}(v_1,v_2)),\quad \Xi_2:=w'\sum_pB_p,\]
\begin{equation}\label{BB}
B_p:=\left\{\begin{array}{ll}
	\hspace{-0.4em}\Theta_2(F^{b'_p}(v'_1,v'_2))K^2(q_r|F^{b_p}(v_1,v_2),F^{b'_p}(v'_1,v'_2)),\quad j\geq q_r,\\
	\hspace{-0.4em} \Theta_2(F^{b'_p}(v'_1,v'_2))K^2(j|F^{b_p}(v_1,v_2),F^{b'_p}(v'_1,v'_2)),\quad j<q_r.
\end{array}\right.
\end{equation}
By (\ref{midex}), we have
\[\Xi_1\leq (1-\exp(-18\Delta(r)))\Delta(r)\sum_p\Theta_2(F^p(v_1,v_2)).\]
By Lemma \ref{keyy1}(2), we write
\begin{equation}\label{lem6i}\Xi_2\leq \Xi_{21}+\Xi_{22},
\end{equation} where
\[\Xi_{21}:= w'\widetilde{\mathbf{K}}_{r}^1(\kappa,R)\sum_{i=0}^{q_R-1}\Theta_2(F^{i}(v'_1,v'_2)),\]
\[\Xi_{22}:= w'\sum_{i=0}^{q_R-1}\Theta_2(F^{i'}(v'_1,v'_2))K^2(F^{i}(v_1,v_2),F^{i'}(v'_1,v'_2)).\]
By Lemma \ref{keyy1}(1) and (\ref{midex}),
\[\Xi_{21}\lesssim w'\varepsilon(\frac{9}{10})^{\frac{\varepsilon R}{2}}\sum_{i=0}^{q_R-1}\Theta_2(F^{i}(v'_1,v'_2))\leq \varepsilon(\frac{9}{10})^{\frac{\varepsilon R}{2}}\sum_{i=0}^{q_R-1}\Theta_2(F^{i}(v_1,v_2)).\]
Note that $e^{\frac{\varepsilon}{10}R}>1$. By Lemma \ref{K22}, Lemma \ref{keyy1}(3),  (\ref{f2x}) and (\ref{midex}),
\[\Xi_{22}\lesssim e^{\frac{\varepsilon}{10}R}q_R^{-\frac{\varepsilon}{3}}w'\sum_{i=0}^{q_R-1}\Theta_2(F^{i'}(v'_1,v'_2))\lesssim (\frac{9}{10})^{\frac{\varepsilon R}{2}}\sum_{i=0}^{q_R-1}\Theta_2(F^{i}(v_1,v_2)).\]
It follows that
\[\Delta(R)-(1-\exp(-18\Delta(r)))\Delta(r)\lesssim (\frac{9}{10})^{\frac{\varepsilon R}{2}},\]
which from the assumption implies $\Delta(R)<\Delta(r)$ and $\Delta(r)\to 0$ as $r\to \infty$ (up to a subsequence). Moreover, we have (\ref{kaka}). This completes the proof of Claim 2.

The estimate (\ref{karrr}) follows from (\ref{kaka}). In fact, we assume (up to a subsequence)
\[\frac{\Delta(r)^2}{(\frac{9}{10})^{\varepsilon r}}\to +\infty.\]
Due to $q_R^{1-\frac{\varepsilon}{2}}\sim q_r$, it follows from Proposition \ref{arthh} that $r<R<\frac{3}{2}r$. Using (\ref{kaka}) repeatedly, we obtain for any $n\in \mathbb{N}$
\[\frac{\Delta(r)^2}{(\frac{9}{10})^{\frac{\varepsilon r}{2}(\frac{3}{4})^n}}\to +\infty,\]
which contradicts the fact that $\Delta(r)\to 0$ as $r\to \infty$. It follows that
\[\Delta(r)\lesssim (\frac{9}{10})^{\frac{\varepsilon r}{2}},\]
which contradicts our assumption at the beginning of the proof.

Using $r<R<\frac{3}{2}r$ again,  we have
\[\Delta(R)\lesssim (\frac{9}{10})^{\varepsilon r}+(\frac{9}{10})^{\frac{\varepsilon R}{2}}\lesssim (\frac{9}{10})^{\frac{\varepsilon R}{2}}.\]
{It follows from (\ref{wan}) that}
\[\widetilde{\mathbf{K}}_R^1(\kappa,R)\lesssim (\frac{9}{10})^{\frac{\varepsilon R}{2}}.\]
This completes the proof of Lemma \ref{karrr1}.
\end{proof}
\subsection{An improvement of the estimate in Lemma \ref{last} }
Based on Lemma \ref{last}, we obtain $\widetilde{\mathbf{K}}_{\bar{N}(\kappa)}^0(\kappa)\lesssim 1$.  It follows that for $0\leq R\leq \tilde{N}(\kappa)$,  the factor $e^{\frac{\varepsilon R}{40}}$ can be replaced by a constant independent of $\kappa$. By repeating the argument in Lemma \ref{keyy1}, we obtain
\[\widetilde{\mathbf{K}}_r^1(\kappa,R)\lesssim \varepsilon q_rq_R^{-1}.\]
Following the same strategy of the proof of Lemma \ref{karrr1}, we have

\begin{Proposition}\label{L1}
	Under the same assumptions of Lemma \ref{last}, we have
	\[\widetilde{\mathbf{K}}_{r}^1(\kappa,r)\lesssim q_r^{-\frac{\varepsilon}{3}},\]
	where $\tilde{N}:=\tilde{N}(\kappa)$ is defined by (\ref{Z12}).
\end{Proposition}

\vspace{2em}

\section{Verification of Criterion 3 for all $\kappa\geq \kappa_0$}\label{S5}
The key ideas of this section are summarized as follows. First, we introduce three  conditions for each $\kappa$:

\begin{itemize}
	\item Condition $R_{\kappa}$:
	\begin{equation}\label{eq:Qk}
		\nabla^1(\kappa,\tilde{N}) \leq \frac{\varepsilon}{960A_\alpha} \quad \text{and} \quad \nabla^2(\kappa,\tilde{N},\tilde{N}) \leq C_0 q_{\tilde{N}}^{-\varepsilon/3}
	\end{equation}
	where $C_0$ is a constant independent of both $\tilde{N}$ and $\varepsilon$.
	\item Condition $S_{\kappa}$:
	\begin{equation}\label{eq:Rk}
		\widetilde{\mathbf{K}}_r^1(\kappa,r) \leq C_1 q_r^{-\varepsilon/3}
	\end{equation}
	where $C_1$ is a constant independent of both $\kappa$ and $r$.
	
	\item Condition $T_{\kappa}$:
	\begin{equation}\label{eq:Tk}
		\Lambda_{II}(\kappa) \leq C_2
	\end{equation}
	where $C_2$ is a constant independent of $\kappa$.
\end{itemize}

In the preceding analysis, we have established the following implication chain:
\begin{equation}\label{f000}
	R_{\kappa} \Rightarrow S_{\kappa} \Rightarrow T_{\kappa}.
\end{equation}

To complete the proof of Theorem~\ref{mathe}, it suffices to verify that $R_{\kappa}$ holds for all $\kappa \geq \kappa_0$. We proceed by induction on $\kappa$:

\begin{enumerate}
	\item  We verify the validity of $R_{\kappa_0}$.
	
	\item  Assuming $R_{\kappa}$ holds for all $\kappa_0 \leq \kappa \leq \kappa_0 + m - 1$, it follows from (\ref{f000}) that both $S_{\kappa}$ and $T_{\kappa}$ hold throughout this range. By using this fact, we obtain that $R_{\kappa_0 + m}$ holds as well.
\end{enumerate}
This induction argument will  complete the proof of our main theorem.

\vspace{1em}

Choose $n_0\in \mathbb{N}$ satisfying the requirement by Lemma \ref{disne}, Lemma \ref{initi},  (\ref{wang1}) and (\ref{W1}) below.
\begin{Lemma}\label{las1}
For each  $\kappa\in \mathbb{N}$ satisfying  $\kappa\geq \kappa_0$, the following hold
\begin{equation}\label{f1}
\nabla^1(\kappa,\tilde{N})\leq \frac{\varepsilon}{960A_\alpha},
\end{equation}
\begin{equation}\label{f2}
\nabla^2(\kappa,\tilde{N},\tilde{N})\leq C_0 q_{\tilde{N}}^{-\frac{\varepsilon}{3}},
\end{equation}
where $\tilde{N}:=\tilde{N}(\kappa)$ is defined by (\ref{Z12}) and $\varepsilon$ is determined by (\ref{epss}).
\end{Lemma}
\begin{Remark}\label{RRX}
	From (\ref{f1}) and (\ref{f2}), the $C^{2+\frac{\varepsilon}{3}}$ smoothness of the preserved invariant circle is evident. For any $\varepsilon' < \varepsilon$, the right-hand side of (\ref{f1}) can be replaced by $\frac{\varepsilon}{k}$, where $k := k(\varepsilon, \varepsilon', \alpha)$ is chosen sufficiently large. By repeating the above proof, the $C^{2+\varepsilon'}$ smoothness of the preserved invariant circle follows directly.
\end{Remark}

We prove Lemma \ref{las1} by induction. At first, we  check the initial step for $\kappa=\kappa_0$. Denote $n_0:=n_{k_0}$ for simplicity and let $v\in (\kappa_0)_{II}$. Let us recall
\[\tilde{N}(k_0)=n_{k_0}-1=n_0-1,\quad s(v)=\frac{r(v)}{\Theta_1(v)}.\] Assume $v$ connects $(x_i,y_i)$ and $(x_j,y_j)$. From Definition \ref{kvec3} and Lemma \ref{disne}, there exists a constant $C_3>0$ such that
\begin{equation}\label{resaa}|r(v)|=|y_j-y_i|\leq C_3q_{n_0}^{-2-\varepsilon}.
\end{equation} By Lemma \ref{initi}, we have $\Theta_1(v)\sim q_{n_0}^{-1}$.
It follows from (\ref{thev1}) that for $(v_1,v_2)\in (\kappa_0, r)_{II}$ and $n_0$ large enough
\begin{equation}\label{wang1}
\nabla^1(v_1,v_2)=\frac{s(v_2)-s(v_1)}{\Theta_2(v_1,v_2)}\leq \frac{s(v_2)-s(v_1)}{\Theta_1(v_1)}\leq C_3 q_{n_0}^{-\varepsilon}\leq \frac{\varepsilon}{1920A_\alpha}.
\end{equation}
For $(v_1,\ldots,v_4)\in (\kappa_0,r,s)_{II}$, we have
\[\lambda_4(v_1,\ldots,v_4)\sim \Theta_4(v_1,\ldots,v_4).\]
Due to the arbitrariness of $v_1$ and $v_2$, we also have for $(v_1,\ldots,v_4)\in (\kappa_0,r,s)_{II}$,
\[\nabla^2(v_1,\ldots,v_4)=\frac{\Theta_4(v_1,\ldots,v_4)}{\lambda_4(v_1,\ldots,v_4)}(\nabla^1(v_3,v_4)-\nabla^1(v_1,v_2))\lesssim q_{n_0}^{-\varepsilon}\leq C_0 q_{\tilde{N}(k_0)}^{-\frac{\varepsilon}{3}}.\]
Note that $\tilde{N}(\kappa)=n_\kappa-2\gamma_0\leq n_\kappa$ for each $\kappa\in \mathbb{N}$. It follows from Lemma \ref{A1} that \[\nabla^1(\kappa_0,\tilde{N}(\kappa_0))\leq \nabla^1(\kappa_0,n_0).\] Hence, (\ref{f1}) and (\ref{f2}) hold for $\kappa=\kappa_0$.

\subsection{Recursive property}
\begin{Lemma}\label{GG}
Let $\Phi(v):=\frac{\Theta_1(v)}{\lambda_1(v)}$. {Given $\kappa \in \mathbb{N}$, if we assume \[
	\nabla^1(\kappa,\tilde{N})\leq \frac{\varepsilon}{960A_\alpha},\quad\nabla^2(\kappa,\tilde{N},\tilde{N})\leq C_0 q_{\tilde{N}}^{-\frac{\varepsilon}{3}}	\]
	Then for each $(v_1,v_2)\in (\kappa,\tilde{N})_{II}$, we have
	\[|\Phi(v_2)-\Phi(v_1)|\leq C_4 q_{\tilde{N}}^{-\frac{\varepsilon}{3}}.\]
	where $C_4$ is independent of $\tilde{N}$ and $\varepsilon$.}
\end{Lemma}
\begin{proof}
By the arithmetics of $\alpha$, one can find an integer $M>0$ independent of $\kappa$ such that
\[\|Mq_{\tilde{N}}\alpha\|-\|q_{\tilde{N}}\alpha\|>\frac{1}{q_{\tilde{N}}}.\]
By the orientation-preserving property (see Proposition \ref{OP}), we have
\[|\Theta_1(v_2)-\Theta_1(v_1)|\leq |\Theta_1(F^{Mq_{\tilde{N}}}(v_1)-\Theta_1(v_1)|.\]
By virtue of Proposition \ref{ko11} and {(\ref{f000})}, for each $v\in (\kappa)_{II}$,
{\[\bold{K}_{\tilde{N}}^0(\kappa)\leq C_1 q_{\tilde{N}}^{-\frac{\varepsilon}{3}}.\]
	It follows that
	\[\exp\left(-C_1q_{\tilde{N}}^{-\frac{\varepsilon}{3}}\right)\leq \frac{\Theta_1(F^{q_{\tilde{N}}}(v))}{\Theta_1(v)}\leq \exp\left(C_1q_{\tilde{N}}^{-\frac{\varepsilon}{3}}\right).\]
	which from the invariance of $\lambda$ with respect to $F$ implies
	\[\exp\left(-MC_1q_{\tilde{N}}^{-\frac{\varepsilon}{3}}\right)\leq \frac{\Theta_1(F^{Mq_{\tilde{N}}}(v_1)}{\Theta_1(v_1)}\leq \exp\left(MC_1q_{\tilde{N}}^{-\frac{\varepsilon}{3}}\right).\]
	It follows that
	\[\frac{|\Theta_1(F^{Mq_{\tilde{N}}}(v_1)-\Theta_1(v_1)|}{\Theta_1(v_1)}\leq \max\left\{\exp\left(MC_1q_{\tilde{N}}^{-\frac{\varepsilon}{3}}\right)-1,1-\exp\left(-MC_1q_{\tilde{N}}^{-\frac{\varepsilon}{3}}\right)\right\},\]
	which gives rise to
	\[|\Theta_1(F^{Mq_{\tilde{N}}}(v_1)-\Theta_1(v_1)|\leq
	2MC_1 q_{\tilde{N}}^{-\frac{\varepsilon}{3}}\Theta_1(v_1).\]
	Hence,
	\[|\Phi(v_2)-\Phi(v_1)|\leq \frac{|\Theta_1(F^{Mq_{\tilde{N}}}(v_1)-\Theta_1(v_1)|}{\lambda(v_1)}\leq 2MC_1 q_{\tilde{N}}^{-\frac{\varepsilon}{3}}\frac{\Theta_1(v_1)}{\lambda(v_1)}\]
	By {(\ref{f000})}, we know that $\frac{\Theta_1(v_1)}{\lambda(v_1)}\leq C_2$. This completes the proof.}
\end{proof}

Based on \cite[Corollary 3.1]{KO1}, we have an improved estimate
\begin{equation}\label{keyfor}
\nabla^2(\kappa,\tilde{N},\tilde{N})\leq C_0q_{\tilde{N}}^{-\eta-\frac{\varepsilon}{3}},
\end{equation}
where $\eta:=\frac{1}{2^{10}}\frac{\varepsilon}{3}$. By using Lemma \ref{GG} and (\ref{keyfor}), we obtain the following recursive property.
\begin{Lemma}\label{SSH}
{Given $\kappa \in \mathbb{N}$, if we assume
	\[|\Phi(v_2)-\Phi(v_1)|\leq C_4 q_{\tilde{N}}^{-\frac{\varepsilon}{3}}.\]
	Then there hold for $l \leq \tilde{N}(\kappa)$,
	\begin{equation}\label{X4}
	\nabla^{1}(\kappa+1,\tilde{N}(\kappa))\leq\nabla^{1}(\kappa,\tilde{N}(\kappa))+10\nabla^{2}(\kappa,\tilde{N}(\kappa),\tilde{N}(\kappa)),
	\end{equation}
	\begin{equation}\label{X1}
	\nabla^{2}(\kappa+1,\tilde{N}(\kappa),l)\leq\nabla^{2}(\kappa,\tilde{N}(\kappa),l)+20\nabla^{2}(\kappa,\tilde{N}(\kappa),\tilde{N}(\kappa)),
	\end{equation}
	\begin{equation}\label{X2}
	\nabla^{1}(\kappa+1,\tilde{N}(\kappa+1))\leq\nabla^{1}(\kappa,\tilde{N}(\kappa))+100(A_{\alpha}+1)^3\nabla^{2}(k,\tilde{N}(\kappa),\tilde{N}(\kappa)).
	\end{equation}}

\end{Lemma}
\begin{proof}
Let $v \in (\kappa+1)_{II}$ and $v_1$ be a contiguous chord in $(\kappa)_{II}$, satisfying $\lambda_1(v) \sim \frac{1}{2}\lambda(v_1)$ as $\kappa\to+\infty$. It follows from the estimate $|\Phi(v_2)-\Phi(v_1)| \leq C_4 q_{\tilde{N}}^{-\varepsilon/3}$ that $\Theta_1(v) \sim \frac{1}{2}\Theta_1(v_1)$.

{There exists  a $\kappa$-chord  $v_2$ such that $(v_1,v_2) \in (\kappa,\tilde{N}(\kappa))_{II}$, which consequently implies $(v+v_1,v_2) \in (\kappa,\tilde{N}(\kappa))_{II}$.} Similar to the proof of Lemma \ref{A1}, we have
\begin{align}
	\nabla^1(v+v_1,v_2) = c\nabla^1(v,v_2) + c_1\nabla^1(v_1,v_2),
\end{align}
where the coefficients are given by:
\[
b =\frac{\Theta_1(v)}{\Theta_1(v)+\Theta_1(v_1)},\quad b_1 = \frac{\Theta_1(v_1)}{\Theta_1(v)+\Theta_1(v_1)}, \]
\[ c = \frac{\Theta_2(v,v_2)}{\Theta_2(v+v_1,v_2)}b, \quad c_1 = \frac{\Theta_2(v_1,v_2)}{\Theta_2(v+v_1,v_2)}b_1.
\]
These satisfy $c + c_1 = 1$ with $c \sim 1/3$ as $\kappa\to+\infty$.

Using the bound
\[
|\nabla^1(v_1,v_2) - \nabla^1(v+v_1,v_2)| \leq \nabla^2(\kappa,\tilde{N}(\kappa),\tilde{N}(\kappa)),
\]
we conclude that $\nabla^1(v,v_2) = c^{-1}\nabla^1(v+v_1,v_2) - (c_1/c)\nabla^1(v_1,v_2)$ remains within $5\nabla^2(\kappa,\tilde{N}(\kappa),\tilde{N}(\kappa))$ of either quantity.

\begin{Remark}
	This establishes that, up to errors bounded by $\nabla^2(\kappa,\tilde{N}(\kappa),\tilde{N}(\kappa))$, the $\nabla^1$ estimates (and similarly for $\nabla^2$) remain valid when replacing one chord in the set $(\kappa)_{II}$  with a chord in the set $(\kappa+1)_{II}$. By iterating this argument, we can successively replace all chords in the set $(\kappa)_{II}$ with chords in the set $(\kappa+1)_{II}$, thereby proving claims (\ref{X4}) and (\ref{X1}).
\end{Remark}

For the claim (\ref{X2}), consider $(v_1,v_2) \in (\kappa+1,\tilde{N}(\kappa+1))_{II}$ and one can find that $v_3 \in (\kappa+1)_{II}$ satisfies $\Theta_2(v_1,v_3) \sim q_{\tilde{N}(\kappa)}^{-1}$. Similar to the proof of Lemma \ref{A1}, we have
\[
\nabla^1(v_1,v_2) = d_1\nabla^1(v_1,v_3) + d_2\nabla^1(v_2,v_3),
\]
with coefficients:
\[
d_1 = -\frac{\Theta_2(v_1,v_3)}{\Theta_2(v_1,v_2)}, \quad d_2 = \frac{\Theta_2(v_2,v_3)}{\Theta_2(v_1,v_2)},
\]
satisfying $d_1 + d_2 = 1$ and $|d_1| \leq q_{\tilde{N}(\kappa+1)}/q_{\tilde{N}(\kappa)}$. Claim (\ref{X2}) then follows from the estimate:
\[
|\nabla^1(v_1,v_3) - \nabla^1(v_2,v_3)| \leq \nabla^2(\kappa+1,\tilde{N}(\kappa),\tilde{N}(\kappa)).
\]
This completes the proof of Lemma \ref{SSH}.		
\end{proof}

\subsection{Proof of Lemma \ref{las1}}
Thanks to (\ref{keyfor}), we have
\[\nabla^2(\kappa,\tilde{N},\tilde{N})\leq C_0q_{\tilde{N}}^{-\eta-\frac{\varepsilon}{3}}.\]
It remains to prove that for each $m>0$,
\begin{equation}\label{Y1}
\nabla^1(\kappa_0+m,\tilde{N}(\kappa_0+m))\leq \frac{\varepsilon}{960A_\alpha}.
\end{equation}
We write \[\tilde{N}_0:=\tilde{N}(\kappa_0),\quad \varepsilon_1:=\frac{\varepsilon}{960A_\alpha},\quad C_5:=100(A_{\alpha}+1)^3C_0.\]
\begin{Lemma}\label{RND}
{For each $ m>1$, if we assume that  for each $\kappa \in [\kappa_0,\kappa_0+m-1]$,
	\[\nabla^1(\kappa,\tilde{N}(\kappa))\leq \frac{\varepsilon}{960A_\alpha},\]}
Then \[{\nabla^1(\kappa_0+m,\tilde{N}(\kappa_0+m))\leq \frac{\varepsilon}{960A\alpha}}.\]
\end{Lemma}

\begin{proof}
{By Lemma \ref{GG} and Lemma \ref{SSH},
	\begin{align*}
		\nabla^1(\kappa_0+1,\tilde{N}(\kappa_0+1))&\leq \nabla^1(\kappa_0,\tilde{N}_0)+100(A_{\alpha}+1)^3\nabla^2(\kappa_0,\tilde{N}_0,\tilde{N}_0)\\
		&\leq \frac{\varepsilon_1}{2}+C_5q_{\tilde{N}_0}^{-\eta-\frac{\varepsilon}{3}}.
	\end{align*}
	Repeating this process, a direct calculation yields
	\begin{equation*}
		\begin{split}
			\nabla^1(\kappa_0+m,\tilde{N}(\kappa_0+m))\leq\frac{\varepsilon_1}{2}+\frac{C_5\log_{2}(2A_\alpha)}{1-\sqrt{2}^{-(\frac{\varepsilon}{3}+\eta)}}\sqrt{2}^{-(\frac{\varepsilon}{3}+\eta)\kappa_0}.
		\end{split}
	\end{equation*}
	Taking the initial data $n_0$ such that
	\begin{equation}\label{W1}
	\kappa_0 > \frac{2\log_{2}\frac{2}{\varepsilon_1}\left(\frac{C_5\log_{2}(2A_\alpha)}{1-\sqrt{2}^{-(\frac{\varepsilon}{3}+\eta)}}\right)}{\frac{\varepsilon}{3}+\eta},
	\end{equation}
	we have
	\[\nabla^1(\kappa_0+m,\tilde{N}(\kappa_0+m))\leq \varepsilon_1=\frac{\varepsilon}{960A_\alpha}.\]
}  \end{proof}
This completes the proof of Lemma \ref{las1}.
\medskip

\appendix
\section{Collection of the notions and notations}
In the order of appearance in the paper, the main notions and notations are listed below for the reader's convenience.

\begin{itemize}
	\item \(A_\alpha\) is the upper bound for the partial quotients in the continued fraction expansion of the constant type rotation number \(\alpha\).
	
	\item \(G_n : \mathbb{T}\times\mathbb{T}\to\mathbb{R}\) stands for the sequence of generating functions constructed in the work.
	
	\item \(F_n : \mathbb{T}\times\mathbb{R}\to\mathbb{T}\times\mathbb{R}\) is the exact area-preserving map generated by \(G_n\). For a fixed \(n\in\mathbb{N}\) we often simply write \(F\) in place of \(F_n\).
	
	\item \(\widetilde{F} : \mathbb{R}^2\to\mathbb{R}^2\) refers to the lift of \(F\).
	
	\item \(\mathcal{O}\) represents a complete orbit of \(F\) with rotation number \(\alpha\).
	
	\item \(\lambda_1(v)\) gives the \(\lambda\)-length of the chord \(v\); \(\Theta_1(v)\) gives the \(\vartheta\)-length of the chord \(v\). Similarly we define \(\lambda_2(v_1,v_2)\), \(\lambda_4(v_1,v_2,v_3,v_4)\) together with \(\Theta_2\), \(\Theta_4\), etc.

	\item \((\kappa)_{\iota}\) designates the set of Type-\(\iota\) \(\kappa\)-chords, where \(\iota\in\{I,II\}\). Accordingly, \((\kappa,r)_{II}\) designates the set of chord pairs \((v_1,v_2)\) belonging to this class, and \((\kappa,r,s)_{II}\) with $s\leq r$ the set of four chords of the same class.
	
	\item $n_\kappa:=\psi(\kappa)$ denotes, in rough terms, the quantity corresponding to $\kappa$ such that $q_{n_\kappa}\sim 2^\kappa$ as $\kappa\to+\infty$; the  definition of $\psi$ is given in (\ref{phius}) and its inverse is denoted by $\phi$.
	
	\item $\tilde{N}(\kappa)=n_\kappa-2\gamma_0$ and $\bar{N}(\kappa)=n_\kappa+2\gamma_0$. These are used to characterize two scales of chord pairs: $(\kappa,\tilde{N}(\kappa))_{II}$ and $(\kappa,\bar{N}(\kappa))_{II}$. A substantial part of the paper is devoted to showing that, under suitable conditions, the estimate on $(\kappa,\tilde{N}(\kappa))_{II}$ can be transferred to $(\kappa,\bar{N}(\kappa))_{II}$. The constant $\gamma_0$, which depends only on $\alpha$; its definition can be found in (\ref{Z11}).
	
	\item  \(\Lambda_\iota(\kappa)\), where \(\iota\in\{I,II\}\), characterise the criteria for the existence of an invariant circle by $F$.

	\item \(\widetilde{\mathbf{K}}_N^0(\kappa)\) measures the supremum of the expansion along orbits generated by \(F\) whose length does not exceed \(q_N\); see the precise definition (\ref{bigK}). When the length equals exactly \(q_N\), we write  \({\mathbf{K}}_N^0(\kappa)\).
	
	\item \(\widetilde{\mathbf{K}}_N^1(\kappa,r)\) measures the supremum of the distortion along orbits generated by \(F\) whose length is at most \(q_N\); the exact definition is given in (\ref{distK}). When the length equals exactly \(q_N\), we write \({\mathbf{K}}_N^1(\kappa,r)\).
	
	\item \(\nabla^1(\kappa,r)\) stands for the supremum of the second-order difference quotient of \(F\) along \(\mathcal{O}\); \(\nabla^2(\kappa,r,s)\) stands for the supremum of the third-order difference quotient of \(F\) along \(\mathcal{O}\).
	\item $C_i$ $(i=0,\ldots,5)$ denotes the positive constant only depending on $\alpha$.
\end{itemize}


\end{document}